\documentclass[abstract,10pt,twocolumn]{scrartcl}

\usepackage[a4paper, left=1.6cm, right=1.6cm, top=2.5cm, bottom=4.0cm]{geometry}
\usepackage[english]{babel}
\usepackage{amsmath}
\usepackage{amssymb}
\usepackage{amsthm}
\usepackage{enumitem}
\usepackage[colorlinks = true, linkcolor = blue, citecolor = blue]{hyperref}
\usepackage[nameinlink]{cleveref}
\usepackage{graphicx}
\usepackage{xcolor}
\usepackage{stfloats}

\setcapindent{0pt}
\setkomafont{caption}{\itshape}
\setkomafont{captionlabel}{\bfseries}

\newtheorem{assumption}{Assumption}
\newtheorem{corollary}{Corollary}
\newtheorem{definizione}{Definition}

\newtheorem{lemma}{Lemma}
\newtheorem{notation}{Notation}
\newtheorem{proposition}{Proposition}
\newtheorem{remark}{Remark}
\newtheorem{theorem}{Theorem}

\crefname{equation}{}{}
\crefname{assumption}{Assumption}{Assumptions}
\crefname{corollary}{Corollary}{Corollaries}
\crefname{definizione}{Definition}{Definitions}
\crefname{example}{Example}{Examples}
\crefname{lemma}{Lemma}{Lemmas}
\crefname{notation}{Notation}{Notation}
\crefname{proposition}{Proposition}{Propositions}
\crefname{remark}{Remark}{Remarks}
\crefname{theorem}{Theorem}{Theorems}

\addto\captionsenglish{}

\title{Extremum Seeking Control for Fully Actuated Mechanical Systems on Lie Groups in the Absence of Dissipation\footnote{This research was supported by the German Research Foundation DFG, project number DA~767/13-1. Corresponding author: Raik~Suttner.}}
\author{Raik Suttner\footnote{Institute of Mathematics, University of Wuerzburg, Wuerzburg, Germany (e-mail: raik.suttner@mathematik.uni-wuerzburg.de)}, $\quad$ Miroslav Krsti\'{c}\footnote{Department of Mechanical and Aerospace Engineering, University of California, San Diego, USA (e-mail: krstic@ucsd.edu)}}
\date{}

\begin{document}
\maketitle
\begin{abstract}
In this paper, we study the problem of extremum seeking control for mechanical systems in dissipation-free environments. This includes attitude control of satellites in space and displacement control of rigid bodies in ideal fluids. The configuration and the velocity of the mechanical system are treated as unknown quantities. The only source of information about the current system state is provided by real-time measurements of a scalar signal whose value has to be minimized. The signal is assumed to be given by a configuration-dependent objective function, which is not known analytically. Our goal is to asymptotically stabilize the mechanical system around states with vanishing velocity and a minimum value of the objective function. The proposed control law employs periodic perturbation signals to extract information about the gradient of the objective function and the velocity of the mechanical system from the response of the sensed signal. Under suitable assumptions, we prove local and non-local stability properties of the closed-loop system. The general results are illustrated by examples.
\end{abstract}

%-----------------------------------------------------------------------------------------------------------------------------------------------------------------
%-----------------------------------------------------------------------------------------------------------------------------------------------------------------
%                                                                Section 1
%-----------------------------------------------------------------------------------------------------------------------------------------------------------------
%-----------------------------------------------------------------------------------------------------------------------------------------------------------------

\section{Introduction}\label{sec:01}
Extremum seeking control for open-loop unstable systems is an important and challenging problem, which has led to extensive research efforts \cite{Zhang20072,Duerr2013,Scheinker2013,Duerr2014,Scheinker2014,Suttner2019}. A popular and frequently studied example is the problem of source seeking with an autonomous agent. In this case, the task is to locate the source of a scalar signal, such as the concentration of a chemical substance or the strength of an electromagnetic field. Depending on the surrounding environment, the autonomous agent can be a robot with wheels~\cite{Zhang20071,Cochran20092}, a drone~\cite{Cochran2009,Matveev2014,Abdelgalil20221}, an underwater vehicle~\cite{Cochran20093,Mandic2020}, or a satellite~\cite{Duerr20142,Walsh2017}. Most of these papers assume that the motion of the agent can be described by a first-order kinematic model. That is, the agent's velocity is controlled directly through the inputs. The same assumption also appears in many studies on extremum seeking control for multi-agent systems~\cite{Stankovic2012,Li2015}. On the other hand, if an agent is controlled through forces or torques, then a second-order dynamic model might be more appropriate than a first-order kinematic model. The method in the present paper is intended for mechanical control systems.

In general, one cannot expect that a mechanical control system is \emph{open-loop stable} in the sense that an arbitrary constant force or torque leads to asymptotic stability. Moreover, because of the system's inertia, vanishing inputs do not necessarily lead to vanishing velocities as in a first-order kinematic model. Thus, the problem of extremum seeking control for second-order dynamic models is more difficult, because it does not only involve the system's configuration (e.g. the position or attitude) but also its velocity. It is clear that stabilizing a mechanical system about an optimal configurations requires some form of dissipation in order to reduce the total energy. Almost all of the existing studies on extremum seeking control for mechanical systems assume that a sufficiently strong loss of energy occurs through velocity-dependent damping~\cite{Michalowsky2014,Michalowsky2015,Scheinker20183,Suttner20192,Suttner20193}. This assumption is justified in environments with significant friction or air resistance. A sufficiently strong damping effect can also be induced through the inputs if measurements of the current velocity are available. However, there are also situations in which none of the above holds; for example, in certain aerospace and underwater applications. The method in the present paper does neither rely on naturally occurring dissipation nor on the availability of velocity measurements.

To the best of our knowledge there are only two studies in the literature so far that address the problem of extremum seeking control for mechanical systems in the absence of naturally occurring dissipation and velocity measurements; namely~\cite{Zhang20072} and~\cite{Suttner20202}. In the present paper, we combine ideas from~\cite{Zhang20072} and~\cite{Suttner20221} to derive a control law for a larger class of second-order dynamic systems. Roughly speaking, our method consists of the following two components: The first component can be seen as a velocity estimator to induce a damping effect and the second component is a gradient estimator to steer the system towards an extremum of the objective function. The idea for the velocity estimator is taken from~\cite{Zhang20072}. The underlying control system in~\cite{Zhang20072} is a double-integrator point mass in the plane. We show that a similar approach as in~\cite{Zhang20072} can be successfully applied to a larger class of fully actuated mechanical systems. Another similarity to the method in~\cite{Zhang20072} is that our control law employs periodic perturbation signals with sufficiently large amplitudes and frequencies to extract gradient information about the objective (gradient estimator). A disadvantage of the method in~\cite{Zhang20072} is that it leads to unbounded velocities in the large-amplitude high-frequency limit. To circumvent this problem, we use a different class of perturbation signals, which ensures bounded velocities. This less invasive perturbation-based approach was already applied in~\cite{Suttner20192} to an acceleration-actuated unicycle, and was extend in~\cite{Suttner20221} to a larger class of mechanical systems. Note, however, that the extremum seeking method in~\cite{Suttner20221} relies on the presence of strict velocity-dependent damping. The approach in~\cite{Suttner20221} does not lead to stability in the absence of naturally occurring dissipation and velocity measurements. We solve the problem in the present paper by combining the gradient estimator from~\cite{Suttner20221} and the velocity estimator from~\cite{Zhang20072} in a suitable way.

As indicated in the previous paragraph, the extremum seeking control law in \cite{Zhang20072} leads to an unbounded growth of velocities with increasing amplitudes and frequencies of the periodic perturbation signals. This undesired feature makes applications to second-order dynamic systems with a non-trivial geometric acceleration impossible. In the present paper, we use a different perturbation-based approach, which is tailored for applications to mechanical systems. Using the averaging theory from~\cite{Bullo2002}, one can show that the closed-loop system approximates the behavior of an averaged system. This in turn leads to the effect that stability properties of the averaged system carry over to the approximating closed-loop system. The same transfer of stability properties also occurs in Lie bracket approximation-based extremum seeking schemes for first-order kinematic systems; see, e.g., \cite{Duerr2013,Duerr2014}. In this case, gradient information is provided by Lie brackets of pairs of suitably chosen vector fields. However, the kinematic approach leads to vanishing Lie brackets if it is applied to second-order dynamic systems. To obtain a suitable method for mechanical systems, we use the class of perturbation signals from~\cite{Bullo2002}. In this case, the averaged system involves so-called \emph{symmetric products} of vector fields from the closed-loop system. One can show that symmetric products originate from iterated Lie brackets of three vector fields on the tangent bundle of the configuration manifold. A suitable design of our extremum seeking method ensures that the symmetric products provide gradient information about the objective function. The symmetric product approach to extremum seeking control is also used in~\cite{Suttner20221} for mechanical systems with strict velocity-dependent dissipation. Here, we go one step further and present a method which does not rely on the presence of dissipation.

The paper is organized as follows. A suitable notion of practical asymptotic stability for the closed-loop system is introduced in \Cref{sec:02}. A precise problem statement, the control law, and the main stability theorems are presented in \Cref{sec:03}. In \Cref{sec:04}, we apply our method to a double-integrator point mass and to a rigid body in an ideal fluid.

%-----------------------------------------------------------------------------------------------------------------------------------------------------------------
%-----------------------------------------------------------------------------------------------------------------------------------------------------------------
%                                                                Section 2
%-----------------------------------------------------------------------------------------------------------------------------------------------------------------
%-----------------------------------------------------------------------------------------------------------------------------------------------------------------

\section{Practical stability}\label{sec:02}
Let $M$ be an embedded submanifold of Euclidean space with Euclidean norm $|\cdot|$. For every $\omega>0$, let $f^\omega$ be a time-dependent vector field on $M$ such that, for every $t_0\in\mathbb{R}$ and every $x_0\in{M}$, the differential equation%
\begin{equation}\label{eq:01}
\dot{x} \ = \ f^\omega(t,x)
\end{equation}%
with initial condition $x(t_0)=x_0$ has a unique maximal solution. The closed-loop system in \Cref{sec:03} will be of the form~\cref{eq:01}, where $\omega$ is a parameter to scale the amplitudes and frequencies of periodic dither signals. For every $\omega>0$ and every $t\in\mathbb{R}$, let $\phi^\omega_{t}\colon{M}\to{M}$ be a bijective map. We use this map to carry out the \emph{change of variables}%
\begin{equation}\label{eq:02}
\tilde{x} \ = \ \phi^{\omega}_t(x).
\end{equation}%
In \Cref{sec:03}, such a change of variables will be applied to the velocities of the closed-loop system.%
\begin{definizione}[\cite{Duerr2013,Duerr2014}]\label{definition1}\em
Let $x_\ast\in{M}$. We say that $x_\ast$ is \emph{practically uniformly stable for~\cref{eq:01} in the variables~\cref{eq:02}} if, for every $\varepsilon>0$, there exist $\omega_0,\delta>0$ such that, for every $\omega\geq\omega_0$, every $t_0\in\mathbb{R}$, and every $\tilde{x}_0\in{M}$ with $|\tilde{x}_0-x_\ast|\leq\delta$, the maximal solution $x$ of~\cref{eq:01} with initial condition $x(t_0)=(\phi^{\omega}_{t_0})^{-1}(\tilde{x}_0)$ satisfies $|\phi^{\omega}_t(x(t))-x_\ast|\leq\varepsilon$ for every $t\geq{t_0}$.
\end{definizione}%
\begin{definizione}[\cite{Duerr2013,Duerr2014}]\label{definition2}\em
Let $x_\ast\in{M}$ and let $S$ be a neighborhood of $x_\ast$ in $M$. We say that $x_\ast$ is \emph{$S$-practically uniformly attractive for~\cref{eq:01} in the variables~\cref{eq:02}} if, for all $r,\varepsilon>0$, there exist $\omega_0,\Delta,R>0$ such that, for every $\omega\geq\omega_0$, every $t_0\in\mathbb{R}$, and every $\tilde{x}_0\in{S}$ with $|\tilde{x}_0-x_\ast|\leq{r}$, the maximal solution $x$ of~\cref{eq:01} with initial condition $x(t_0)=(\phi^{\omega}_{t_0})^{-1}(\tilde{x}_0)$ satisfies $|\phi^{\omega}_t(x(t))-x_\ast|\leq{R}$ for every $t\geq{t_0}$ and $|\phi^{\omega}_t(x(t))-x_\ast|\leq\varepsilon$ for every $t\geq{t_0+\Delta}$.
\end{definizione}%
\begin{definizione}[\cite{Duerr2013,Duerr2014}]\label{definition3}\em
Let $x_\ast\in{M}$ and let $S$ be a neighborhood of $x_\ast$ in $M$. We say that $x_\ast$ is \emph{$S$-practically uniformly asymptotically stable for~\cref{eq:01} in the variables~\cref{eq:02}} if $x_\ast$ is practically uniformly stable for~\cref{eq:01} in the variables~\cref{eq:02} and if $x_\ast$ is $S$-practically uniformly attractive for~\cref{eq:01} in the variables~\cref{eq:02}.
\end{definizione}%
\begin{remark}\label{remark1}\em
If the conditions in \Cref{definition2,definition3} are satisfied for some unknown (possibly small) neighborhood $S$ of $x_\ast$ in $M$, then we replace the prefix ``$S$-'' by the word ``locally''. If there exists a time-dependent vector field $\bar{f}$ on $M$ such that $f^\omega=\bar{f}$ for every $\omega>0$, then we omit the word ``practically'' in \Cref{definition1,definition2,definition3}. If $\phi^{\omega}_t(x)=x$ for every $\omega>0$, every $t\in\mathbb{R}$, and every $x\in{M}$, then we omit the phrase ``in the variables~\cref{eq:02}'' in \Cref{definition1,definition2,definition3}. The word ``uniformly'' in \Cref{definition1,definition2,definition3} indicates that the properties therein are uniform with respect to the time parameter.
\end{remark}%
Let $\bar{f}$ be a time-dependent vector field on $M$ such that, for every $t_0\in\mathbb{R}$ and every $\bar{x}_0\in{M}$, the differential equation%
\begin{equation}\label{eq:03}
\dot{\bar{x}} \ = \ \bar{f}(t,\bar{x})
\end{equation}%
with initial condition $\bar{x}(t_0)=\bar{x}_0$ has a unique maximal solution. In \Cref{sec:03}, the averaged system of the closed-loop system will be of the form~\cref{eq:03}.%
\begin{definizione}[\cite{Duerr2013,Duerr2014}]\label{definition4}\em
We say that \emph{the solutions of~\cref{eq:01} in the variables~\cref{eq:02} approximate the solutions of~\cref{eq:03}} if, for every compact subset $K$ of $M$ and all $\nu,\Delta>0$, there exists $\omega_0>0$ such that, for every $t_0\in\mathbb{R}$ and every $\bar{x}_0\in{K}$, the following implication holds: If the maximal solution $\bar{x}$ of~\cref{eq:03} with initial condition $\bar{x}(t_0)=\bar{x}_0$ satisfies $\bar{x}(t)\in{K}$ for every $t\in[t_0,t_0+\Delta]$, then, for every $\omega\geq\omega_0$, the maximal solution $x$ of~\cref{eq:01} with initial condition $x(t_0)=(\phi^{\omega}_{t_0})^{-1}(\bar{x}_0)$ satisfies $|\phi^{\omega}_t(x(t))-\bar{x}(t)|\leq\nu$ for every $t\in[t_0,t_0+\Delta]$.
\end{definizione}%
We will see in \Cref{sec:03} that the closed-loop system approximated the behavior of an averaged system in the sense of \Cref{definition4}. In the situation of \Cref{definition4}, stability properties of~\cref{eq:03} carry over to the approximating system~\cref{eq:01} as follows.%
\begin{proposition}[\cite{Duerr2013,Duerr2014}]\label{proposition1}
Suppose that the embedded submanifold $M$ is a topologically closed subset of the ambient Euclidean space.\footnote{Erratum: The assumption that $M$ is a topologically closed subset of Euclidean space is missing in Proposition~1 in \cite{Suttner20221}.} Assume that the solutions of~\cref{eq:01} in the variables~\cref{eq:02} approximate the solutions of~\cref{eq:03}. Let $x_\ast\in{M}$ and let $S$ be a neighborhood of $x_\ast$ in $M$. If $x_\ast$ is $S$-uniformly asymptotically stable for~\cref{eq:03}, then $x_\ast$ is $S$-practically uniformly asymptotically stable for~\cref{eq:01} in the variables~\cref{eq:02}.
\end{proposition}%
A proof of \Cref{proposition1} can be found in \Cref{sec:appendixA}.

In \Cref{sec:03}, we will apply \Cref{proposition1} to conclude stability properties of the closed-loop system from stability properties of the averaged system.

%-----------------------------------------------------------------------------------------------------------------------------------------------------------------
%-----------------------------------------------------------------------------------------------------------------------------------------------------------------
%                                                                Section 3
%-----------------------------------------------------------------------------------------------------------------------------------------------------------------
%-----------------------------------------------------------------------------------------------------------------------------------------------------------------

\section{Problem statement, control law, and main results}\label{sec:03}
Throughout this section, we suppose that%
\begin{itemize}
	\item $G$ is a Lie group of dimension $n$ such that $G$ is a closed embedded submanifold of Euclidean space,
	\item $\mathbb{I}$ is an inner product on the tangent space~$\mathfrak{g}$ to $G$ at the identity,
	\item $R$ is a symmetric and positive-semidefinite endomorphism of~$\mathfrak{g}$ with respect to~$\mathbb{I}$,
	\item $\psi$ is a smooth real-valued function on $G$.
\end{itemize}%
We use a similar notation as in the textbook~\cite{BulloBook}.%
\begin{notation}\label{notation1}\em
For every $g\in{G}$, let $L_g$ denote the \emph{left translation map by $g$} and let $T_eL_g$ denote the restriction of its tangent map to $\mathfrak{g}$. For every $v\in\mathfrak{g}$, let $v^{\text{L}}$ denote the left-invariant vector field on $G$ defined by $v^{\text{L}}(g):=T_eL_g(v)$. Let $\nabla$ denote the left-invariant \emph{Levi-Civita connection} for the left-invariant Riemannian metric on $G$ induced by $\mathbb{I}$. The \emph{restriction of $\nabla$ to $\mathfrak{g}$} is the well-defined bilinear map $\stackrel{\mathfrak{g}}{\nabla}\colon\mathfrak{g}\times\mathfrak{g}\to\mathfrak{g}$, written $(v,w)\mapsto\stackrel{\mathfrak{g}}{\nabla}_v\!w$, such that%
\begin{equation}\label{eq:04}
T_eL_g\big(\stackrel{\mathfrak{g}}{\nabla}_v\!w\big) \ = \ \big(\nabla_{v^{\text{L}}}w^{\text{L}}\big)(g)
\end{equation}%
for every $g\in{G}$ and all $v,w\in\mathfrak{g}$; cf.~Theorem 5.40 in~\cite{BulloBook}.
\end{notation}%

%-----------------------------------------------------------------------------------------------------------------------------------------------------------------
%                                                                Section 3.1
%-----------------------------------------------------------------------------------------------------------------------------------------------------------------

\subsection{Problem statement}\label{sec:03:01}
We assume that we deal with a fully actuated control system on the \emph{configuration manifold} $G$ whose dynamics can be described by \emph{Euler--Poincar{\'e} equations} of the form%
\begin{equation}\label{eq:05}
\dot{g} \ = \ T_eL_g(v), \qquad \dot{v} + \stackrel{\mathfrak{g}}{\nabla}_v\!v \ = \ -Rv + \mathrm{u}
\end{equation}%
on $G\times\mathfrak{g}$, where we use \Cref{notation1}. The first equation in~\cref{eq:05} describes the \emph{velocity} $\dot{g}$ in the current \emph{configuration}~$g$. Since $G$ is a Lie group, we can represent $\dot{g}$ equivalently through the \emph{body velocity} $v:=(T_eL_g)^{-1}(\dot{g})$, which is $\mathfrak{g}$-valued. The expression on the left-hand side of the second equation in~\cref{eq:05} represents the \emph{geometric acceleration of $g$ with respect to~$\nabla$} in terms of the body velocity. We assume that the geometric acceleration can be controlled directly through a $\mathfrak{g}$-valued \emph{input}~$\mathrm{u}$. \emph{Velocity-dependent dissipation} is described by the symmetric and negative-semidefinite endomorphism~$-R$ on $\mathfrak{g}$. Note that, in particular, we allow vanishing $R$, which means absence of dissipation.

The configuration~$g$ and the velocity~$v$ in \cref{eq:05} as well as the map $R$ are unknown quantities. The only information about the current system state is provided by real-time measurements of the \emph{output}%
\begin{equation}\label{eq:06}
\mathrm{y} \ = \ \psi(g).
\end{equation}%
Note that a measurement of~\cref{eq:06} does not provide immediate information about the velocity, because the \emph{objective function}~$\psi$ only depends on the configuration. We are interested in an output feedback control law~$\mathrm{u}$ for~\cref{eq:05} that stabilizes the closed-loop system around states $(g_\ast,v_\ast)$, where $g_\ast\in{G}$ is a minimizer of $\psi$ and $v_\ast\in\mathfrak{g}$ is zero. Reaching zero velocity becomes especially difficult in the absence of dissipative forces (i.e. for vanishing~$R$). A suitable control law is proposed in the next subsection.

%-----------------------------------------------------------------------------------------------------------------------------------------------------------------
%                                                                Section 3.2
%-----------------------------------------------------------------------------------------------------------------------------------------------------------------

\subsection{Control law}\label{sec:03:02}
Before we present our control law for~\cref{eq:05}, we explain the underlying idea and indicate how it works. For the sake of simplicity, we restrict our introductory discussion to the particular example of a double-integrator point in the Euclidean space $\mathbb{R}^n$ whose dynamics are described by the second-order differential equation%
\begin{equation}\label{eq:07}
\ddot{g} \ = \ -r\,\dot{g} + \mathrm{u},
\end{equation}%
where $r$ is some nonnegative damping constant. The input $\mathrm{u}$ is $\mathbb{R}^n$-valued and the output is given by~\cref{eq:06}. The first component of our control law is the extremum seeking method from~\cite{Suttner20221}. This method is an output feedback strategy, which employs perturbation signals with sufficiently large amplitudes and high frequencies. If we apply the control law from~\cite{Suttner20221} to the above double-integrator, then the solutions of the closed-loop system approximate the solutions of an averaged system in the sense of \Cref{definition4}. By \Cref{proposition1}, this implies that stability properties of the averaged system carry over to the approximating closed-loop system. It is shown in~\cite{Suttner20221} that the averaged system is driven into the direction of the negative gradient of $\psi$ up to some varying scaling factor. If we omit this scaling factor for the sake of simplicity, then the averaged system is of the form%
\begin{equation}\label{eq:08}
\ddot{\bar{g}} \ = \ -r\,\dot{\bar{g}} - \lambda\,\mathrm{grad}\psi(\bar{g}),
\end{equation}%
where $\lambda$ is a positive constant and $\mathrm{grad}\psi$ denotes the gradient of $\psi$ with respect to the standard Euclidean inner product on $\mathbb{R}^n$. If $r>0$, then the averaged system is a damped oscillator with ``potential function'' $\psi$ and, under mild assumptions on $\psi$, one can show that minimizers of $\psi$ lead to asymptotically stable equilibria. The situation changes for $r=0$; i.e., in the absence of dissipation. Then, the averaged system is an undamped oscillator and we cannot expect asymptotic stability. Since we do not assume measurements of the velocity, we do not have direct access to the current negative velocity $-\dot{\bar{g}}(t)$ for the purpose of linear damping. To circumvent this problem, we use an approach from~\cite{Zhang20072}, which induces damping as follows. For a solution $t\mapsto\bar{g}(t)$ of the oscillator \cref{eq:08}, we can expect that the signal $t\mapsto-\dot{\bar{g}}(t)$ leads the signal $t\mapsto-\mathrm{grad}\psi(\bar{g}(t))$ with a certain phase shift. Thus, if we want access to $t\mapsto-\dot{\bar{g}}(t)$, then we should add phase to $t\mapsto-\mathrm{grad}\psi(\bar{g}(t))$. Now, instead of~\cref{eq:08}, we consider the system%
\begin{equation}\label{eq:09}
\ddot{\bar{g}} \ = \ - r\,\dot{\bar{g}} + \mathrm{y}_{\text{c}},
\end{equation}%
where $\mathrm{y}_{\text{c}}$ is the output of the phase-lead compensator%
\begin{subequations}\label{eq:10}%
\begin{align}
\dot{\bar{w}} & \ = \ - a\,\bar{w} + \kappa\,\mathrm{u}_{\text{c}}, \allowdisplaybreaks \\
\mathrm{y}_{\text{c}} & \ = \ - b\,\bar{w} + \lambda\,\mathrm{u}_{\text{c}}
\end{align}%
\end{subequations}%
with input $\mathrm{u}_{\text{c}}=-\mathrm{grad}\psi(\bar{g})$ and positive gains $a$, $\lambda$, $b$, $\kappa$ such that%
\begin{equation}\label{eq:11}
a\,\lambda - b\,\kappa \ > \ 0.
\end{equation}%
The transfer function $G_{\text{c}}$ of phase-lead compensator the~\cref{eq:10} is given by $G_{\text{c}}(z)=\frac{\lambda{z}+a\lambda-b\kappa}{z+a}$. We combine \cref{eq:09} and \cref{eq:10} and obtain the system%
\begin{subequations}\label{eq:12}%
\begin{align}
\ddot{\bar{g}} & \ = \ -r\,\dot{\bar{g}} - b\,\bar{w} - \lambda\,\mathrm{grad}\psi(\bar{g}), \label{eq:09:a} \allowdisplaybreaks \\
\dot{\bar{w}} & \ = \ -a\,\bar{w} - \kappa\,\mathrm{grad}\psi(\bar{g}). \label{eq:09:b}
\end{align}%
\end{subequations}%
Under mild assumptions on $\psi$, one can show that minimizers of $\psi$ lead to asymptotically stable equilibria for~\cref{eq:12} for every $r\geq0$; in particular for $r=0$. Thus, if we can approximate the solutions of an averaged system of the form \cref{eq:12} by means of a suitable control law, then we can expect practical stability for the closed-loop system. This is exactly what we do in the next paragraph; not only for a double-integrator point but for the more general control system~\cref{eq:05}.

Now we introduce the components of the control law. Let $T$ be a positive real number and let $\boldsymbol{u}$ be a measurable and bounded $\mathbb{R}^n$-valued map on $\mathbb{R}$. As in~\cite{Suttner20221}, the component functions $u^1,\ldots,u^n$ of $\boldsymbol{u}$ shall act as independent periodic perturbation signals. For this reason, we demand that $\boldsymbol{u}$ is $T$-periodic and zero-mean, and that the component functions $U^1,\ldots,U^n$ of the zero-mean antiderivative $\boldsymbol{U}$ of $\boldsymbol{u}$ satisfy the orthonormality condition%
\begin{equation}\label{eq:13}
\int_0^TU^i(\tau)\,U^j(\tau)\,\mathrm{d}\tau \ = \ \left\{ \begin{tabular}{cl} $1$ & if $i=j$, \\ $0$ & if $i\neq{j}$ \end{tabular} \right.
\end{equation}%
for all $i,j\in\{1,\ldots,n\}$. For example, for $T=2\pi$, we can define the component functions of $\boldsymbol{u}$ by%
\begin{equation}\label{eq:14}
u^i(\tau) \ := \ \sqrt{2}\,i\,\cos(i\tau)
\end{equation}%
for every $i\in\{1,\ldots,n\}$. Next, let $\alpha$ be a smooth real-valued function on $\mathbb{R}$. As in~\cite{Suttner20221}, the function $\alpha$ is a design function for the output feedback, which must have the following properties. We demand that the product of $\alpha$ and its derivative $\alpha'$ is strictly increasing and attains only positive values. For example, we can define $\alpha$ by%
\begin{equation}\label{eq:15}
\alpha(z) \ := \ \sqrt{z+\log(2\cosh{z})},
\end{equation}%
because then $\alpha\alpha'=(1+\tanh)/2$ has the desired properties. Let $a$, $\lambda$, $b$, $\kappa$ be positive real numbers such that \cref{eq:11} is satisfied. Choose an orthonormal basis $(e_1,\ldots,e_n)$ for~$\mathfrak{g}$ with respect to~$\mathbb{I}$. Recall that the restriction $\stackrel{\mathfrak{g}}{\nabla}$ of the Levi-Civita connection to $\mathfrak{g}$ is introduced in \Cref{notation1}. Set%
\begin{equation}\label{eq:16}
\mu \ := \ \stackrel{\mathfrak{g}}{\nabla}_{e_1}\!e_1 + \cdots + \stackrel{\mathfrak{g}}{\nabla}_{e_n}\!e_n \in \mathfrak{g}.
\end{equation}%
To relate elements of $\mathfrak{g}$ to elements of $\mathbb{R}^n$, we introduce the following convenient notation.%
\begin{notation}\label{notation2}\em
Let $(\boldsymbol{e}_1,\ldots,\boldsymbol{e}_n)$ denote the standard basis for $\mathbb{R}^n$. Let $E\colon\mathbb{R}^n\to\mathfrak{g}$ denote the vector space isomorphism that assigns to each basis vector $\boldsymbol{e}_i$ of $\mathbb{R}^n$ the corresponding basis vector $e_i$ of $\mathfrak{g}$. We use the following notation to represent the values of $E$ and its inverse $E^{-1}$. If the boldface symbol $\boldsymbol{w}$ is a vector in $\mathbb{R}^n$, then we denote $E(\boldsymbol{w})$ by the lightface symbol $w$. Conversely, if the lightface symbol $w$ is a vector in $\mathfrak{g}$, then we denote $E^{-1}(w)$ by the boldface symbol $\boldsymbol{w}$. For any such pair of $\boldsymbol{w}\in\mathbb{R}^n$ and $w\in\mathfrak{g}$, we can write $\boldsymbol{w}=w^1\boldsymbol{e}_1+\cdots+w^n\boldsymbol{e}_n$ and $w=w^1e_1+\cdots+w^ne_n$ with the same components $w^1,\ldots,w^n\in\mathbb{R}$. In the same way, if $\boldsymbol{\mathrm{u}}$ is an $\mathbb{R}^n$-valued map, then we denote the $\mathfrak{g}$-valued composition $E\circ\boldsymbol{\mathrm{u}}$ by the symbol $\mathrm{u}$ and vice versa. 
\end{notation}%
Using \Cref{notation2}, we propose the parameter- and time-dependent output feedback control law%
\begin{subequations}\label{eq:17}%
\begin{align}
\boldsymbol{\mathrm{u}} & \ = \ - b\,(\boldsymbol{w}-\kappa\,\alpha(\mathrm{y}-\eta)\,\boldsymbol{U}(\omega{t})) \label{eq:17:a} \allowdisplaybreaks \\
& \qquad + \lambda\,\alpha(\mathrm{y}-\eta)\,\omega\,\boldsymbol{u}(\omega{t}) + \lambda^2\,\alpha^2(\mathrm{y}-\eta)\,\boldsymbol{\mu} \label{eq:17:b}
\end{align}%
\end{subequations}%
for~\cref{eq:05} with output $\mathrm{y}$ as in~\cref{eq:06}, where $\omega$ is a positive real parameter, $\boldsymbol{w}$ is the $\mathbb{R}^n$-valued state of the ``phase-lead compensator''%
\begin{subequations}\label{eq:18}%
\begin{align}
\dot{\boldsymbol{w}} & \ = \ -a\,(\boldsymbol{w}-\kappa\,\alpha(\mathrm{y}-\eta)\,\boldsymbol{U}(\omega{t})) \label{eq:18:a} \allowdisplaybreaks \\
& \qquad + \kappa\,\alpha(\mathrm{y}-\eta)\,\omega\,\boldsymbol{u}(\omega{t}) \label{eq:18:b}
\end{align}%
\end{subequations}%
to induce damping in the averaged system, and $\eta$ is the real-valued state of the high-pass filter%
\begin{equation}\label{eq:19}
\dot{\eta} \ = \ -h\,\eta + h\,\mathrm{y}
\end{equation}%
with filter output $\mathrm{y}-\eta$ and positive gain $h$ to remove a possible offset from $\mathrm{y}$. Now we apply \cref{eq:17,eq:18,eq:19} to~\cref{eq:05} and, using \Cref{notation2}, we obtain the closed-loop system%
\begin{subequations}\label{eq:20}%
\begin{align}
\dot{g} & \ = \ T_eL_g(v), \label{eq:20:a} \allowdisplaybreaks \\
\!\!\dot{v} + \stackrel{\mathfrak{g}}{\nabla}_v\!v & \ = \ -Rv - b\,(w-\kappa\,\alpha(\psi(g)-\eta)\,U(\omega{t})) \label{eq:20:b} \allowdisplaybreaks \\
& \qquad + \lambda\,\alpha(\psi(g)-\eta)\,\omega\,u(\omega{t}) \label{eq:20:c} \allowdisplaybreaks \\
& \qquad + \lambda^2\,\alpha^2(\psi(g)-\eta)\,\mu, \label{eq:20:d} \allowdisplaybreaks \\
\dot{\boldsymbol{w}} & \ = \ -a\,(\boldsymbol{w}-\kappa\,\alpha(\mathrm{y}-\eta)\,\boldsymbol{U}(\omega{t})), \label{eq:20:e} \allowdisplaybreaks \\
& \qquad + \kappa\,\alpha(\psi(g)-\eta)\,\omega\,\boldsymbol{u}(\omega{t}), \label{eq:20:f} \allowdisplaybreaks \\
\dot{\eta} & \ = \ -h\,\eta + h\,\psi(g) \label{eq:20:g}
\end{align}%
\end{subequations}%
on the state manifold $M:=G\times\mathfrak{g}\times\mathbb{R}^n\times\mathbb{R}$. To present the averaged system of~\cref{eq:20}, we introduce the following notation. Let $\mathrm{grad}\psi$ denote the gradient vector field of $\psi$ with respect to the left-invariant Riemannian metric on $G$ induced by~$\mathbb{I}$. Using \Cref{notation1}, we define a $\mathfrak{g}$-valued map $\mathfrak{grad}\psi$ on $G$ by%
\begin{equation}\label{eq:21}
\mathfrak{grad}\psi(\bar{g}) \ := \ (T_eL_g)^{-1}(\mathrm{grad}\psi(\bar{g})).
\end{equation}%
Using \Cref{notation2}, for every $\bar{g}\in{G}$, we denote the representation of $\mathfrak{grad}\psi(\bar{g})$ as an element of $\mathbb{R}^n$ by $\boldsymbol{\mathfrak{grad}}\psi(\bar{g})$. If we carry out the change of variables%
\begin{subequations}\label{eq:22}%
\begin{align}
\tilde{g} & \ = \ g, \label{eq:22:a} \allowdisplaybreaks \\
\tilde{v} & \ = \ v - \lambda\,\alpha(\psi(g)-\eta)\,U(\omega{t}), \label{eq:22:b} \allowdisplaybreaks \\
\tilde{\boldsymbol{w}} & \ = \ \boldsymbol{w} - \kappa\,\alpha(\psi(g)-\eta)\,\boldsymbol{U}(\omega{t}), \label{eq:22:c} \allowdisplaybreaks \\
\tilde{\eta} & \ = \ \eta \label{eq:22:d}
\end{align}%
\end{subequations}%
for the state of \cref{eq:20}, then a standard first-order averaging argument can be applied to extract the averaged system%
\begin{subequations}\label{eq:23}%
\begin{align}
\dot{\bar{g}} & \ = \ T_eL_g(\bar{v}), \label{eq:23:a} \allowdisplaybreaks \\
\dot{\bar{v}} + \stackrel{\mathfrak{g}}{\nabla}_{\bar{v}}\!\bar{v} & \ = \ -R\bar{v} - b\,\bar{w} \label{eq:23:b} \allowdisplaybreaks \\
& \qquad - \lambda^2\,(\alpha\alpha')(\psi(\bar{g})-\bar{\eta})\,\mathfrak{grad}\psi(\bar{g}), \label{eq:23:c} \allowdisplaybreaks \\
& \!\!\!\!\!\!\!\!\!\!\!\!\!\!\!\!\!\!\!\! \dot{\bar{\boldsymbol{w}}} \ = \ -a\,\bar{\boldsymbol{w}} - \kappa\,\lambda\,(\alpha\alpha')(\psi(\bar{g})-\bar{\eta})\,\boldsymbol{\mathfrak{grad}}\psi(\bar{g}), \label{eq:23:d} \allowdisplaybreaks \\
\dot{\bar{\eta}} & \ = \ -h\,\bar{\eta} + h\,\psi(\bar{g}) \label{eq:23:e}
\end{align}%
\end{subequations}%
on $M$. Alternatively, one can also apply the more general averaging theory from~\cite{Bullo2002} to come to the same conclusion. In the terminology of \Cref{definition4}, we can state the following approximation result.%
\begin{proposition}\label{proposition2}
The solutions of~\cref{eq:20} in the variables \cref{eq:22} approximate the solutions of~\cref{eq:23}.
\end{proposition}%
A proof of \Cref{proposition2} can be found in \Cref{sec:appendixB}.%
\begin{remark}\label{remark2}\em
Using the averaging procedure from~\cite{Bullo2002}, one can show that the gradient term in~\cref{eq:23:c} originates from so-called symmetric products of vector fields which are given by the term $\alpha(\psi(g)-\eta)$ in \cref{eq:20:c}. To make this statement more precise, recall (e.g. from \cite{BulloBook}) that, for all smooth vector fields $X,Y$ on $G$, the \emph{symmetric product of $X,Y$ with respect to $\nabla$} is the vector field%
\begin{equation}\label{eq:24}
\langle{X}\colon\!{Y}\rangle \ := \ \nabla_XY + \nabla_YX
\end{equation}%
on $G$. One can show that the oscillatory term in~\cref{eq:20:c} leads to the following sum of symmetric products in the averaged system:%
\begin{subequations}\label{eq:25}%
\begin{align}
& -\frac{1}{2}\,\sum_{i=1}^n\langle\lambda\,\alpha(\psi-\eta)\,e_i^{\text{L}}\colon\!\lambda\,\alpha(\psi-\eta)\,e_i^{\text{L}}\rangle \label{eq:25:a} \allowdisplaybreaks \\
&  \ = \ -\lambda^2\,(\alpha\alpha')(\psi-\eta)\,\mathrm{grad}\psi - \lambda^2\,\alpha^2(\psi-\eta)\,\mu^{\text{L}}, \label{eq:25:b}
\end{align}%
\end{subequations}%
where we use \Cref{notation1}. The first term in \cref{eq:25:b} corresponds to the one in \cref{eq:23:c}. The second term in \cref{eq:25:b} can be seen as an undesired contribution, which however does not appear in~\cref{eq:23} since it is compensated by~\cref{eq:20:d}. Note that the state manifold of a mechanical system is the tangent bundle of its configuration manifold. It is known from~\cite{Bullo2002} that the symmetric products in the averaged system originate from certain iterated Lie brackets on the tangent bundle of the configuration manifold. For this reason, one may say that the proposed extremum seeking method is based on approximations of Lie brackets. This differential geometric interpretation of control law~\cref{eq:17,eq:18,eq:19} establishes a direct link to the approach in \cite{Duerr2013,Duerr2014} for kinematic control systems.
\end{remark}%
Because of \Cref{proposition1,proposition2}, stability properties of the averaged system~\cref{eq:23} carry over to the approximating closed-loop system~\cref{eq:20}. We will use this property in the next subsection to derive the main stability results for the closed-loop system.

%-----------------------------------------------------------------------------------------------------------------------------------------------------------------
%                                                                Section 3.3
%-----------------------------------------------------------------------------------------------------------------------------------------------------------------

\subsection{Main results}\label{sec:03:03}
It is clear that a stability result for the closed-loop system~\cref{eq:20} requires certain assumptions on the objective function~$\psi$. A local stability property can be proved if $\psi$ is locally quadratic around a minimizer. This ensures that the negative gradient of $\psi$ provides a sufficiently strong force towards a minimum.%
\begin{assumption}\label{assumption1}
The objective function $\psi$ attains a local minimum value at some point $g_\ast$ of $G$ and the Hessian of $\psi$ at $g_\ast$ is positive definite.
\end{assumption}%
Under the above assumption, we can show that the averaged system~\cref{eq:23} is locally asymptotically stable, and therefore, by \Cref{proposition1,proposition2}, the following practical stability property holds for the closed-loop system on the state manifold $M=G\times\mathfrak{g}\times\mathbb{R}^n\times\mathbb{R}$ in the terminology of \Cref{definition3,remark1}.%
\begin{theorem}\label{theorem1}
Suppose that \Cref{assumption1} is satisfied with $g_\ast$ as therein. Then the point $(g_\ast,0,\boldsymbol{0},\psi(g_\ast))$ of $M$ is locally practically uniformly asymptotically stable for~\cref{eq:20} in the variables~\cref{eq:22}.
\end{theorem}%
A proof of \Cref{theorem1} is given in \Cref{sec:appendixC}.

Next, we extend \Cref{theorem1} to a non-local stability result and provide an estimate for the domain of attraction. This, however, requires that the underlying control system is ``well-behaved'' in the sense that its geometric acceleration does \emph{not} display quadratic growth with increasing velocities. Note that, for large velocities, the bilinear term on the left-hand side of~\cref{eq:23:b} may have a larger magnitude than the negative gradient force on the right-hand side of~\cref{eq:23:c}. In this case, one cannot expect that the estimator in~\cref{eq:23:d} provides reliable information about the current body velocity. This in turn can lead to a loss of the desired damping effect. For this reason, we assume the following property in our non-local stability result (see \Cref{theorem2} below).%
\begin{assumption}\label{assumption2}
We have $\stackrel{\mathfrak{g}}{\nabla}_v\!v=0$ for every $v\in\mathfrak{g}$.
\end{assumption}%
In some texts, e.g.~in~\cite{PostnikovBook}, \Cref{assumption2} is the defining property of a so-called \emph{Cartan connection}.%
\begin{remark}[Lemma~7.2 in~\cite{Milnor1976}]\label{remark3}\em
Suppose that $G$ is connected. Then \Cref{assumption2} is satisfied if and only if the left-invariant Riemannian metric on $G$ is actually bi-invariant.
\end{remark}%
If \Cref{assumption2} is satisfied, then we can also weaken the requirement of a positive definite Hessian in \Cref{assumption1} to the condition that the gradient of $\psi$ is nonvanishing around a minimizer. To quantify the set with a nonvanishing gradient, we introduce the following notation. For every manifold $X$, every real-valued function $f$ on $X$, every point $x$ of $X$, and every real number $y>f(x)$, let $f^{-1}(\leq{y},x)$ denote the connected component of the $y$-sublevel set of $f$ containing $x$.%
\begin{assumption}\label{assumption3}
There exist $g_\ast\in{G}$ and $y_0>\psi(g_\ast)$ such that $\psi^{-1}(\leq{y_0},g_\ast)$ is compact and $\mathrm{grad}\psi(g)\neq0$ for every $g\in\psi^{-1}(\leq{y_0},g_\ast)$ with $g\neq{g_\ast}$.
\end{assumption}%
To provide an estimate for the domain of attraction of the closed-loop system, we introduce the following maps. For every symmetric and positive-semidefinite endomorphism $\Xi$ of $\mathfrak{g}$ with respect to $\mathbb{I}$, define a real-valued function $\|\cdot\|_\Xi$ on $\mathfrak{g}$ by $\|v\|_\Xi:=\sqrt{\mathbb{I}(\Xi{v},v)}$. In particular, for the identity map $I$ of $\mathfrak{g}$, we get the norm $\|\cdot\|:=\|\cdot\|_{I}$ induced by $\mathbb{I}$. Because of inequality~\cref{eq:11},%
\begin{equation}\label{eq:26}
\Xi \ := \ b\,\kappa\,\big((a\,\lambda-b\,\kappa)\,I + \lambda\,R\big)^{-1}
\end{equation}%
\begin{figure*}%
\centering$\begin{matrix}\includegraphics{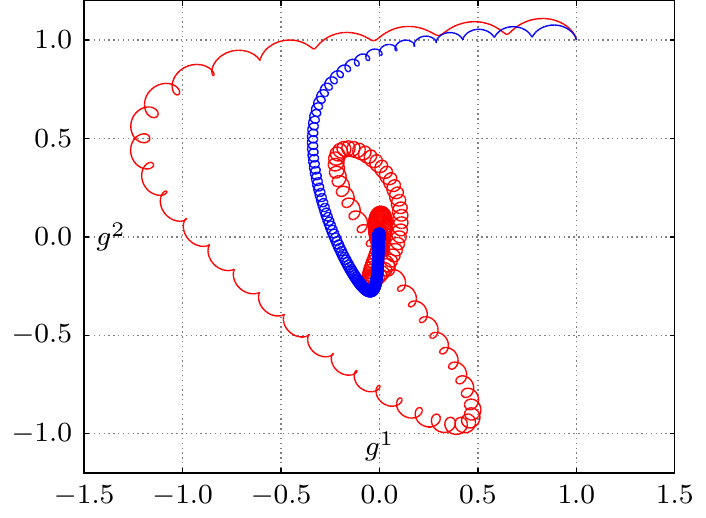}\end{matrix}\qquad\qquad\begin{matrix}\includegraphics{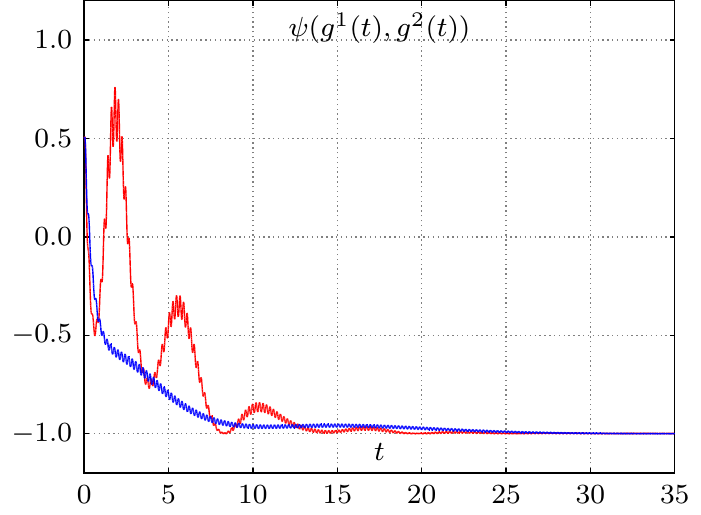}\end{matrix}$%
\caption{Position (left) and output (right) of a source-seeking double-integrator point in the plane. The simulation results are generated with a ``fortunate'' (blue) and a ``realistic'' (red) choice of the parameters $a$, $\lambda$, $b$, $\kappa$. The ``fortunate'' choice of $a$, $\lambda$, $b$, $\kappa$ results in a monotonic decay of the output without overshoots of the source.}%
\label{figure1}%
\end{figure*}%
is a well-defined symmetric and positive-definite endomorphism of $\mathfrak{g}$ with respect to $\mathbb{I}$. Define a real-valued function $\beta$ on $\mathbb{R}$ by%
\begin{equation}\label{eq:27}
\beta(z) \ := \ \int_0^z\big((\alpha\alpha')(\tilde{z})-(\alpha\alpha')(0)\big)\,\mathrm{d}\tilde{z}.
\end{equation}%
The assumptions on $\alpha$ in \Cref{sec:03:02} imply that $\beta$ is positive definite about~$0$. Finally, if $g_\ast$ is a minimum point of $\psi$, define a real-valued function $V$ on $M$ by%
\begin{subequations}\label{eq:28}%
\begin{align}
& V(g,v,\boldsymbol{w},\eta) \ := \ \tfrac{1}{2}\,\|v\|^2 + \tfrac{1}{2}\,\big\|v-\tfrac{\lambda}{\kappa}\,w\big\|_{\Xi}^2 \label{eq:28:a} \allowdisplaybreaks \\
& + \lambda^2\,(\alpha\alpha')(0)\,(\psi(g)-\psi(g_\ast)) + \lambda^2\,\beta(\psi(g)-\eta). \label{eq:28:b}
\end{align}%
\end{subequations}%
The terms on the right-hand sides of \cref{eq:28:a} and~\cref{eq:28:b} can be interpreted as \emph{kinetic energy} and \emph{potential energy}, respectively.\footnote{Erratum: The term ``$-\psi(g_\ast)$'' is missing in the definitions of the Lyapunov functions in~\cite{Suttner20221}.} Since the system state $w$ induces damping, one may refer to the second term in \cref{eq:28:a} as a \emph{damping-modified kinetic energy}. Using the sublevel sets of~$V$, we can state the following non-local stability result for the closed-loop system~\cref{eq:20} on the state manifold $M=G\times\mathfrak{g}\times\mathbb{R}^n\times\mathbb{R}$ in the terminology of \Cref{definition3}.%
\begin{theorem}\label{theorem2}
Suppose that \Cref{assumption2,assumption3} are satisfied with $g_\ast$ and $y_0$ as therein. Set $x_\ast:=(g_\ast,0,\boldsymbol{0},\psi(g_\ast))$ and $V_0:=\lambda^2(\alpha\alpha')(0)(y_0-\psi(g_\ast)) $. Then the point $x_\ast$ of $M$ is $V^{-1}(\leq{V_0},x_\ast)$-practically uniformly asymptotically stable for~\cref{eq:20} in the variables~\cref{eq:22}.
\end{theorem}%
A proof of \Cref{theorem2} is given in \Cref{sec:appendixC}.

%-----------------------------------------------------------------------------------------------------------------------------------------------------------------
%-----------------------------------------------------------------------------------------------------------------------------------------------------------------
%                                                                Section 4
%-----------------------------------------------------------------------------------------------------------------------------------------------------------------
%-----------------------------------------------------------------------------------------------------------------------------------------------------------------

\section{Examples}\label{sec:04}

%-----------------------------------------------------------------------------------------------------------------------------------------------------------------
%                                                                Section 4.1
%-----------------------------------------------------------------------------------------------------------------------------------------------------------------

\subsection{Double integrator point in the plane}\label{sec:04:01}
As an easy example, we consider a double-integrator point on $G=\mathbb{R}^2$ in the absence of dissipative forces. This problem is also studied in~\cite{Zhang20072}. To allow a direct comparison with the control law in~\cite{Zhang20072}, we choose the perturbation signals $u^1:=-\sqrt{2}\sin$ and $u^2:=\sqrt{2}\cos$. Then the zero-mean antiderivatives $U^1:=\sqrt{2}\cos$ and $U^2:=\sqrt{2}\sin$ satisfy the orthonormality condition~\cref{eq:13}. Choose an orthonormal basis $(e_1,e_2)$ for $\mathfrak{g}=\mathbb{R}^2$ with respect to the Euclidean inner product $\mathbb{I}$. Then, we can write elements $g\in{G}$ and $v\in\mathfrak{g}$ component-wise as $g=g^1e_1+g^2e_2$ and $v=v^1e_2+v^2e_2$, respectively. In the same way, we can write a vector $\boldsymbol{w}$ in $\mathbb{R}^2$ component-wise as $\boldsymbol{w}=w^1\boldsymbol{e}_1+w^2\boldsymbol{e}_2$ with respect to the standard basis $(\boldsymbol{e}_1,\boldsymbol{e}_2)$ for $\mathbb{R}^2$. For the double-integrator point, the general closed-loop system~\cref{eq:20} reduces to%
\begin{subequations}\label{eq:29}%
\begin{align}
\dot{g}^1 & \ = \ v^1, \qquad \dot{g}^2 \ = \ v^2, \label{eq:29:a} \allowdisplaybreaks \\
\dot{v}^1 & \ = \ - b\,(w^1-\kappa\,\alpha(\psi(g)-\eta)\,\sqrt{2}\cos(\omega{t})) \label{eq:29:b} \allowdisplaybreaks \\
& \qquad - \lambda\,\alpha(\psi(g)-\eta)\,\omega\,\sqrt{2}\sin(\omega{t}), \label{eq:29:c} \allowdisplaybreaks \\
\dot{v}^2 & \ = \ - b\,(w^2-\kappa\,\alpha(\psi(g)-\eta)\,\sqrt{2}\sin(\omega{t})) \label{eq:29:d} \allowdisplaybreaks \\
& \qquad + \lambda\,\alpha(\psi(g)-\eta)\,\omega\,\sqrt{2}\cos(\omega{t}), \label{eq:29:e} \allowdisplaybreaks \\
\dot{w}^1 & \ = \ - a\,(w^1-\kappa\,\alpha(\psi(g)-\eta)\,\sqrt{2}\cos(\omega{t})) \label{eq:29:f} \allowdisplaybreaks \\
& \qquad - \kappa\,\alpha(\psi(g)-\eta)\,\omega\,\sqrt{2}\sin(\omega{t}), \label{eq:29:g} \allowdisplaybreaks \\
\dot{w}^2 & \ = \ - a\,(w^2-\kappa\,\alpha(\psi(g)-\eta)\,\sqrt{2}\sin(\omega{t})) \label{eq:29:h} \allowdisplaybreaks \\
& \qquad + \kappa\,\alpha(\psi(g)-\eta)\,\omega\,\sqrt{2}\cos(\omega{t}), \label{eq:29:i} \allowdisplaybreaks \\
\dot{\eta} & \ = \ -h\,\eta + h\,\psi(g) \label{eq:29:j}
\end{align}%
\end{subequations}%
and the general change of variables~\cref{eq:22} reduces to%
\begin{subequations}\label{eq:30}%
\begin{align}
\tilde{g}^1 & \ = \ g^1, \qquad \tilde{g}^2 \ = \ g^2, \qquad \tilde{\eta} \ = \ \eta, \label{eq:30:a} \allowdisplaybreaks \\
\tilde{v}^1 & \ = \ v^1 - \lambda\,\alpha(\psi(g)-\eta)\,\sqrt{2}\cos(\omega{t}), \label{eq:30:b} \allowdisplaybreaks \\
\tilde{v}^2 & \ = \ v^2 - \lambda\,\alpha(\psi(g)-\eta)\,\sqrt{2}\sin(\omega{t}), \label{eq:30:c} \allowdisplaybreaks \\
\tilde{w}^1 & \ = \ w^1 - \kappa\,\alpha(\psi(g)-\eta)\,\sqrt{2}\cos(\omega{t}), \label{eq:30:d} \allowdisplaybreaks \\
\tilde{w}^2 & \ = \ w^2 - \kappa\,\alpha(\psi(g)-\eta)\,\sqrt{2}\sin(\omega{t}). \label{eq:30:e}
\end{align}%
\end{subequations}%
Note that, in this particular example, \Cref{assumption2} is trivially satisfied and the energy-like function $V$ defined in~\cref{eq:28} is given by%
\begin{subequations}\label{eq:31}%
\begin{align}
& V(g,v,\boldsymbol{w},\eta) = \frac{1}{2}\,|v|^2 + \frac{b\,\kappa}{2\,(a\,\lambda-b\,\kappa)}\,\big|v-\tfrac{\lambda}{\kappa}\,w\big|^2 \label{eq:31:a} \allowdisplaybreaks \\
& + \lambda^2\,(\alpha\alpha')(0)\,(\psi(g) - \psi(g_\ast)) + \lambda^2\,\beta(\psi(g)-\eta), \label{eq:31:b}
\end{align}%
\end{subequations}%
where $|\cdot|$ denotes the Euclidean norm on $\mathbb{R}^2$. As a particular case of \Cref{theorem2}, we get the following result.%
\begin{corollary}\label{corollary1}
Suppose that \Cref{assumption3} is satisfied with $g_\ast$ and $y_0$ as therein. Set $x_\ast:=(g_\ast,0,\boldsymbol{0},\psi(g_\ast))$ and $V_0:=\lambda^2(\alpha\alpha')(0)(y_0 - \psi(g_\ast))$. Then the point $x_\ast$ of $\mathbb{R}^2\times\mathbb{R}^2\times\mathbb{R}^2\times\mathbb{R}$ is $V^{-1}(\leq{V_0},x_\ast)$-practically uniformly asymptotically stable for~\cref{eq:29} in the variables~\cref{eq:30}.
\end{corollary}%
To allow a visual comparison with the simulation results in~\cite{Zhang20072}, we test our method for the same quadratic objective function as therein. That is, we define $\psi$ by%
\begin{equation}\label{eq:32}
\psi(g^1,g^2) \ := \ -1 + (g^1)^2 + (g^2)^2/2.
\end{equation}%
In this case, \Cref{corollary1} implies that $x_\ast=(0,0,\boldsymbol{0},-1)$ is semi-globally practically asymptotically stable for the closed-loop system. We define $\alpha$ by \cref{eq:15} and, as in~\cite{Zhang20072}, we set $\omega:=30$ and $h:=1$. In the first simulation, we choose $a:=\lambda:=b:=1$ and $\kappa:=1/2$ to satisfy inequality~\cref{eq:11}. The position and the output of the double-integrator point are shown in red in \Cref{figure1}. We observe that the position converges into a neighborhood of the origin. We also see that the trajectory overshoots the optimal point several times with oscillations of period $\approx5$. For a second simulation, we determine the linearization of the averaged system about $x_\ast$ and compute the eigenvalues as a function of the parameters $a$, $\lambda$, $b$, $\kappa$. It turns out that the imaginary parts of the eigenvalues are very small for $a:=1.3$, $\lambda:=0.7$, $b:=1.2$, and $\kappa:=0.7$. The corresponding results are shown in blue in \Cref{figure1}. For this choice of $a$, $\lambda$, $b$, $\kappa$, there are almost no swings and overshoots. Note, however, that this selection of $a$, $\lambda$, $b$, $\kappa$ requires knowledge about $\psi$.

%-----------------------------------------------------------------------------------------------------------------------------------------------------------------
%                                                                Section 4.2
%-----------------------------------------------------------------------------------------------------------------------------------------------------------------

\subsection{Rigid body in an ideal fluid}\label{sec:04:02}
In our second example, we consider a system on the special Euclidean group $G=\mathrm{SE}(3)=\mathrm{SO}(3)\ltimes\mathbb{R}^3$. Then, $G$ is a $6$-dimensional embedded submanifold of the $12$-dimensional Euclidean space $\mathbb{R}^{3\times{3}}\times\mathbb{R}^3$ and its tangent space $\mathfrak{g}$ at the identity is $\mathfrak{so}(3)\times\mathbb{R}^3$, where $\mathfrak{so}(3)$ denotes the subspace of skew-symmetric matrices in $\mathbb{R}^{3\times{3}}$. We refer the reader to Example~5.49 in the textbook~\cite{BulloBook} for more details. It is convenient to identify elements of $\mathfrak{so}(3)$ with elements of $\mathbb{R}^3$ through the vector space isomorphism%
\begin{equation}\label{eq:33}
\widehat{\cdot}\colon\mathbb{R}^3\to\mathfrak{so}(3), \qquad \left[\begin{smallmatrix} \Omega^1 \\ \Omega^2 \\ \Omega^3 \end{smallmatrix}\right] \mapsto \left[\begin{smallmatrix} 0 & -\Omega^3 & \Omega^2 \\ \Omega^3 & 0 & -\Omega^1 \\ -\Omega^2 & \Omega^1 & 0 \end{smallmatrix}\right].
\end{equation}%
This in turn induces the vector space isomorphism $(\boldsymbol{\Omega},\boldsymbol{V})\mapsto(\widehat{\boldsymbol{\Omega}},\boldsymbol{V})$ between $\mathbb{R}^3\oplus\mathbb{R}^3$ and $\mathfrak{g}$. From now on we represent elements of $\mathfrak{g}$ by elements of $\mathbb{R}^3\oplus\mathbb{R}^3$, where elements in $\mathbb{R}^3$ are written as column vectors. Let $\mathbb{J}$ and $\mathbb{M}$ be two inner products on $\mathbb{R}^3$. Let $[\mathbb{J}]$ and $[\mathbb{M}]$ denote the matrix representations of $\mathbb{J}$ and $\mathbb{M}$ with respect to the standard basis of $\mathbb{R}^3$. Typical elements of $G=\mathrm{SO}(3)\ltimes\mathbb{R}^3$ and $\mathfrak{g}\cong\mathbb{R}^3\oplus\mathbb{R}^3$ are denoted by $(\boldsymbol{R},\boldsymbol{r})$ and $(\boldsymbol{\Omega},\boldsymbol{V})$, respectively. Let $\times$ denote the usual cross product on $\mathbb{R}^3$. We consider the \emph{controlled Kirchhoff equations}%
\begin{subequations}\label{eq:34}%
\begin{align}
\dot{\boldsymbol{R}} \ = \ \boldsymbol{R}\,\widehat{\boldsymbol{\Omega}}, \qquad \dot{\boldsymbol{r}} & \ = \ \boldsymbol{R}\,\boldsymbol{V}, \label{eq:34:a} \allowdisplaybreaks \\
\dot{\boldsymbol{\Omega}} + [\mathbb{J}]^{-1}(\boldsymbol{\Omega}\times[\mathbb{J}]\boldsymbol{\Omega} + \boldsymbol{V}\times[\mathbb{M}]\boldsymbol{V}) & \ = \ [\mathbb{J}]^{-1}\boldsymbol{\mathrm{f}}_{\boldsymbol{\Omega}}, \label{eq:34:b} \allowdisplaybreaks \\
\dot{\boldsymbol{V}} + [\mathbb{M}]^{-1}(\boldsymbol{\Omega}\times[\mathbb{M}]\boldsymbol{V}) & \ = \ [\mathbb{M}]^{-1}\boldsymbol{\mathrm{f}}_{\boldsymbol{V}}, \label{eq:34:c}
\end{align}%
\end{subequations}%
where $\boldsymbol{\mathrm{f}}_{\boldsymbol{\Omega}}$ and $\boldsymbol{\mathrm{f}}_{\boldsymbol{V}}$ are $\mathbb{R}^3$-valued inputs. The above Kirchhoff equations are also used in other works, such as \cite{FossenBook,Leonard1997,Bullo1999}, as a simple model for the motion of a rigid body in an ideal fluid. Here, the term \emph{ideal fluid} means that there are no dissipative forces. Let $\check{\cdot}\colon\mathfrak{so}(3)\to\mathbb{R}^3$ denote the inverse of $\hat{\cdot}$. The Kirchhoff equations represent Euler--Poincar{\'e} equations of the form~\cref{eq:05} with $R=0$ if we carry the linear input transformation%
\begin{equation}\label{eq:35}
\boldsymbol{\mathrm{f}}_{\boldsymbol{\Omega}} \ = \ [\mathbb{J}]\check{\mathrm{u}}_\Omega  \qquad \text{and} \qquad \boldsymbol{\mathrm{f}}_{\boldsymbol{V}} \ = \ [\mathbb{M}]\mathrm{u}_V
\end{equation}%
with a $\mathfrak{g}$-valued input $\mathrm{u}=(\mathrm{u}_\Omega,\mathrm{u}_V)$ consisting of an $\mathfrak{so}(3)$-valued component $\mathrm{u}_\Omega$ and an $\mathbb{R}^3$-valued component $\mathrm{u}_V$. An explicit formula for the restriction $\stackrel{\mathfrak{g}}{\nabla}$ of the Levi-Civita connection to $\mathfrak{g}$ is given by equation~(5.27) in~\cite{BulloBook}.

Let $\psi$ be a smooth real-valued function on $\mathrm{SE}(3)$. We want to apply the output feedback control law \cref{eq:17}. To this end, we choose an $\mathbb{R}^6$-valued map~$\boldsymbol{u}$ of periodic perturbation signals with zero-mean antiderivatives $\boldsymbol{U}$ and a smooth real-valued function $\alpha$ as in \Cref{sec:03:02}. We decompose $\boldsymbol{u}$ and $\boldsymbol{U}$ as $\boldsymbol{u}=(\boldsymbol{u}_{\boldsymbol{\Omega}},\boldsymbol{u}_{\boldsymbol{V}})$ and $\boldsymbol{U}=(\boldsymbol{U}_{\boldsymbol{\Omega}},\boldsymbol{U}_{\boldsymbol{V}})$ with $\mathbb{R}^3$-valued maps $\boldsymbol{u}_{\boldsymbol{\Omega}}$, $\boldsymbol{u}_{\boldsymbol{V}}$, $\boldsymbol{U}_{\boldsymbol{\Omega}}$, and $\boldsymbol{U}_{\boldsymbol{V}}$. Let $a$, $\lambda$, $b$, $\kappa$ be positive real-numbers such that~\cref{eq:11} holds. Following the construction in \Cref{sec:03:02}, we have to choose an orthonormal basis for $\mathfrak{g}$. For this purpose, we let $[\mathbb{J}]^{-1/2}$ and $[\mathbb{M}]^{-1/2}$ denote the inverses of the square roots of the symmetric and positive-definite matrices $[\mathbb{J}]$ and $[\mathbb{M}]$, respectively. Let $(\boldsymbol{e}_1,\boldsymbol{e}_2,\boldsymbol{e}_3)$ denote the standard basis for $\mathbb{R}^3$. For every $i\in\{1,2,3\}$, set $e_i:=(\widehat{[\mathbb{J}]^{-1/2}\boldsymbol{e}_i},\boldsymbol{0}_{3\times{1}})$ and $e_{i+3}:=(\boldsymbol{0}_{3\times{3}},[\mathbb{M}]^{-1/2}\boldsymbol{e}_i)$. Then $(e_1,\ldots,e_6)$ is an orthonormal basis for $\mathfrak{g}$. A direct computation, using equation~(5.27) in~\cite{BulloBook}, shows that the element $\mu$ of $\mathfrak{g}$ in~\cref{eq:16} is zero. Let $E\colon\mathbb{R}^6\to\mathfrak{g}$ be the isomorphism from \Cref{notation2} with respect to $(e_1,\ldots,e_6)$. We take the control law $\boldsymbol{\mathrm{u}}$ in \cref{eq:17} and apply $\mathrm{u}=E\circ\boldsymbol{u}$ after the linear transformation~\cref{eq:35} to~\cref{eq:34}. For the rigid body in an ideal fluid, the general closed-loop system~\cref{eq:20} reduces to%
\begin{figure*}%
\centering$\begin{matrix}\begin{matrix}\includegraphics{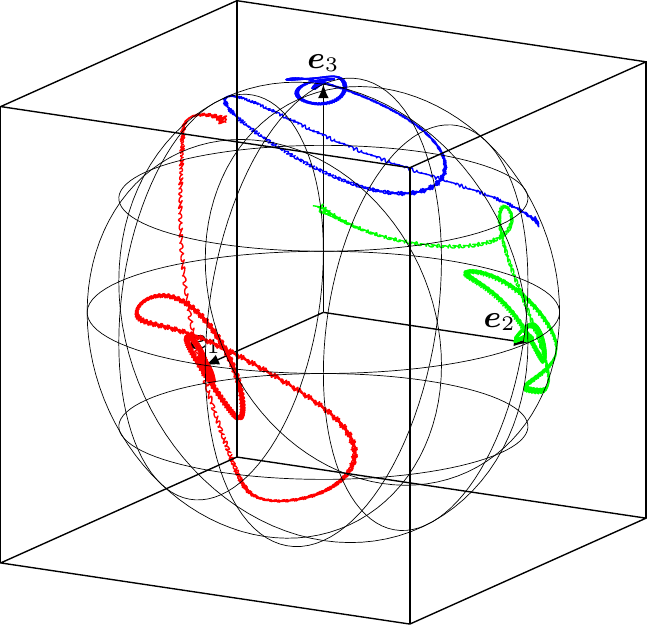}\end{matrix} & \qquad\qquad & \begin{matrix}\includegraphics{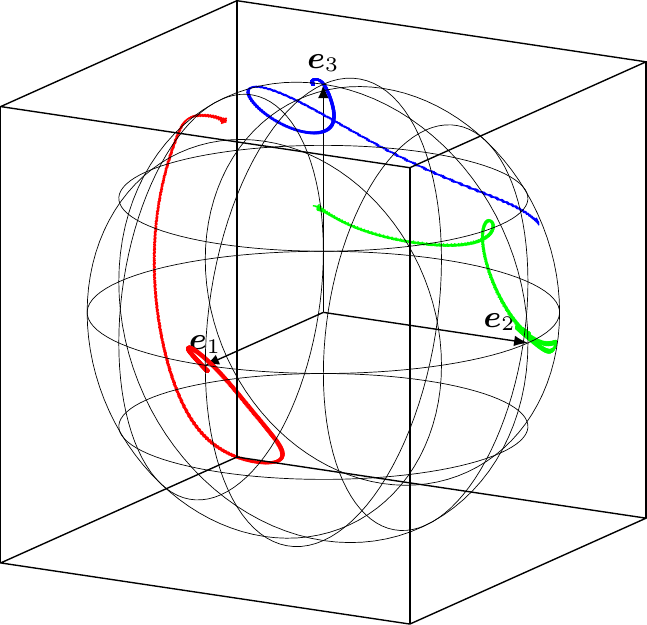}\end{matrix}\\\begin{matrix}\includegraphics{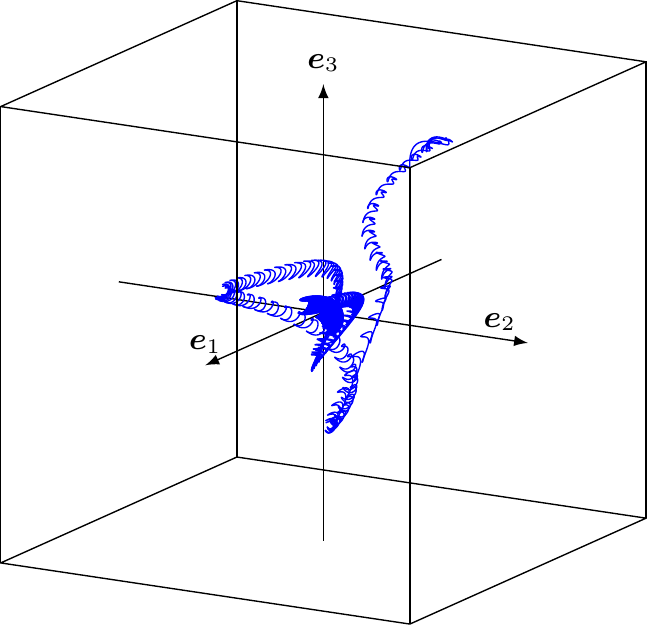}\end{matrix} & \qquad\qquad & \begin{matrix}\includegraphics{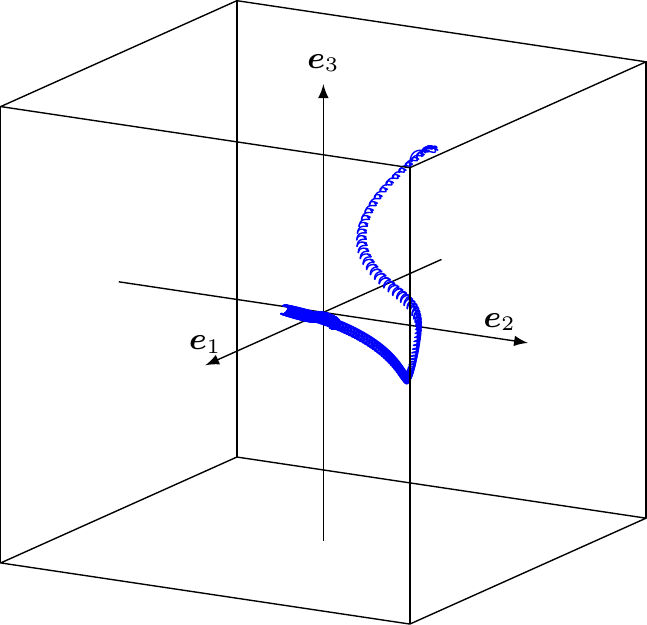}\end{matrix}\end{matrix}$%
\caption{Attitude (upper row) and position (lower row) of a rigid body in an ideal fluid. The results in the left column are generated for $a:=b:=\kappa:=0.1$ and $\lambda:=0.2$ and the results in the right column are generated for $a:=0.15$ and $\lambda:=b:=\kappa:=0.13$. The latter choice of $a$, $\lambda$, $b$, $\kappa$ leads to smaller imaginary parts of the eigenvalues of the linearization of the averaged system and therefore to smaller oscillations and overshoots of the optimal state. The attitudes are illustrated by the curves $t\mapsto\boldsymbol{R}(t)\boldsymbol{e}_1$ (red), $t\mapsto\boldsymbol{R}(t)\boldsymbol{e}_2$ (green), $t\mapsto\boldsymbol{R}(t)\boldsymbol{e}_3$ (blue).}%
\label{figure2}%
\end{figure*}%
\begin{subequations}\label{eq:36}%
\begin{align}
\dot{\boldsymbol{R}} & \ = \ \boldsymbol{R}\,\widehat{\boldsymbol{\Omega}}, \qquad \dot{\boldsymbol{r}} \ = \ \boldsymbol{R}\,\boldsymbol{V}, \label{eq:36:a} \allowdisplaybreaks \\
\dot{\boldsymbol{\Omega}} & + [\mathbb{J}]^{-1}(\boldsymbol{\Omega}\times[\mathbb{J}]\boldsymbol{\Omega} + \boldsymbol{V}\times[\mathbb{M}]\boldsymbol{V}) \label{eq:36:b} \allowdisplaybreaks \\
& \ = \ + \lambda\,\alpha(\psi(\boldsymbol{R},\boldsymbol{r})-\eta)\,\omega\,[\mathbb{J}]^{-1/2}\boldsymbol{u}_{\boldsymbol{\Omega}}(\omega{t}) \label{eq:36:c} \allowdisplaybreaks \\
- b\,& [\mathbb{J}]^{-1/2}\big(\boldsymbol{w}_{\boldsymbol{\Omega}}-\kappa\,\alpha(\psi(\boldsymbol{R},\boldsymbol{r})-\eta)\,\boldsymbol{U}_{\boldsymbol{\Omega}}(\omega{t})\big), \label{eq:36:d} \allowdisplaybreaks \\
\dot{\boldsymbol{V}} & + [\mathbb{M}]^{-1}(\boldsymbol{\Omega}\times[\mathbb{M}]\boldsymbol{V}) \label{eq:36:e} \allowdisplaybreaks \\
& \ = \ + \lambda\,\alpha(\psi(\boldsymbol{R},\boldsymbol{r})-\eta)\,\omega\,[\mathbb{M}]^{-1/2}\boldsymbol{u}_{\boldsymbol{V}}(\omega{t}) \label{eq:36:f} \allowdisplaybreaks \\
- b\,& [\mathbb{M}]^{-1/2}\big(\boldsymbol{w}_{\boldsymbol{V}}-\kappa\,\alpha(\psi(\boldsymbol{R},\boldsymbol{r})-\eta)\,\boldsymbol{U}_{\boldsymbol{V}}(\omega{t})\big), \label{eq:36:g} \allowdisplaybreaks \\
\dot{\boldsymbol{w}}_{\boldsymbol{\Omega}} & \ = \ - a\,\big(\boldsymbol{w}_{\boldsymbol{\Omega}}-\kappa\,\alpha(\psi(\boldsymbol{R},\boldsymbol{r})-\eta)\,\boldsymbol{U}_{\boldsymbol{\Omega}}(\omega{t})\big) \label{eq:36:h} \allowdisplaybreaks \\
& \qquad + \kappa\,\alpha(\psi(\boldsymbol{R},\boldsymbol{r})-\eta)\,\omega\,\boldsymbol{u}_{\boldsymbol{\Omega}}(\omega{t}), \label{eq:36:i} \allowdisplaybreaks \\
\dot{\boldsymbol{w}}_{\boldsymbol{V}} & \ = \ - a\,\big(\boldsymbol{w}_{\boldsymbol{V}}-\kappa\,\alpha(\psi(\boldsymbol{R},\boldsymbol{r})-\eta)\,\boldsymbol{U}_{\boldsymbol{V}}(\omega{t})\big) \label{eq:36:j} \allowdisplaybreaks \\
& \qquad + \kappa\,\alpha(\psi(\boldsymbol{R},\boldsymbol{r})-\eta)\,\omega\,\boldsymbol{u}_{\boldsymbol{V}}(\omega{t}), \label{eq:36:k} \allowdisplaybreaks \\
\dot{\eta} & \ = \ -h\,\eta + h\,\psi(\boldsymbol{R},\boldsymbol{r}) \nonumber \tag{36l}
\end{align}%
\end{subequations}%
and the general change of variables~\cref{eq:22} reduces to%
\begin{subequations}\label{eq:37}%
\begin{align}
\tilde{\boldsymbol{R}} & \ = \ \boldsymbol{R}, \qquad \tilde{\boldsymbol{r}} \ = \ \boldsymbol{r}, \qquad \tilde{\eta} \ = \ \eta, \label{eq:37:a} \allowdisplaybreaks \\
\tilde{\boldsymbol{\Omega}} & \ = \ \boldsymbol{\Omega} - \lambda\,\alpha(\psi(\boldsymbol{R},\boldsymbol{r})-\eta)\,[\mathbb{J}]^{-1/2}\boldsymbol{U}_{\boldsymbol{\Omega}}(\omega{t}), \label{eq:37:b} \allowdisplaybreaks \\
\tilde{\boldsymbol{V}} & \ = \ \boldsymbol{V} - \lambda\,\alpha(\psi(\boldsymbol{R},\boldsymbol{r})-\eta)\,[\mathbb{M}]^{-1/2}\boldsymbol{U}_{\boldsymbol{V}}(\omega{t}), \label{eq:37:c} \allowdisplaybreaks \\
\tilde{\boldsymbol{w}}_{\boldsymbol{\Omega}} & \ = \ \boldsymbol{w}_{\boldsymbol{\Omega}} - \kappa\,\alpha(\psi(\boldsymbol{R},\boldsymbol{r})-\eta)\,\boldsymbol{U}_{\boldsymbol{\Omega}}(\omega{t}), \label{eq:37:d} \allowdisplaybreaks \\
\tilde{\boldsymbol{w}}_{\boldsymbol{V}} & \ = \ \boldsymbol{w}_{\boldsymbol{V}} - \kappa\,\alpha(\psi(\boldsymbol{R},\boldsymbol{r})-\eta)\,\boldsymbol{U}_{\boldsymbol{V}}(\omega{t}). \label{eq:37:e}
\end{align}%
\end{subequations}%
As a particular case of \Cref{theorem1}, we get the following stability result for the closed-loop system~\cref{eq:36} on the state manifold $M=\mathrm{SE}(3)\times\mathbb{R}^6\times\mathbb{R}^6\times\mathbb{R}$.%
\begin{corollary}\label{corollary2}
Suppose that \Cref{assumption1} is satisfied with $g_\ast=(\boldsymbol{R}_\ast,\boldsymbol{r}_\ast)\in\mathrm{SE}(3)$ as therein. Then the point $(g_\ast,\boldsymbol{0},\boldsymbol{0},\psi(g_\ast))$ of $M$ is locally practically uniformly asymptotically stable for~\cref{eq:36} in the variables~\cref{eq:37}.
\end{corollary}%
We test our method in the following situation. The inner products $\mathbb{J}$ and $\mathbb{M}$ are defined through their representation matrices with respect to the standard basis by%
\begin{equation}\label{eq:38}
[\mathbb{J}] := \tfrac{1}{3}\left[\begin{smallmatrix} 5 & 0 & -2 \\ 0 & 7 & 2 \\ -2 & 2 & 6 \end{smallmatrix}\right] \quad \text{and} \quad [\mathbb{M}] := \tfrac{1}{3}\left[\begin{smallmatrix} 7 & 0 & 2 \\ 0 & 5 & -2 \\ 2 & -2 & 6 \end{smallmatrix}\right].
\end{equation}%
The objective function $\psi$ on $\mathrm{SE}(3)$ is defined by%
\begin{equation}\label{eq:39}
\psi(\boldsymbol{R},\boldsymbol{r}) \ := \ |\boldsymbol{R}-\boldsymbol{I}_3|^2/4 + |\boldsymbol{r}|^2/2,
\end{equation}%
where $\boldsymbol{I}_3\in\mathrm{SO}(3)$ denotes the $3\times{3}$ identity matrix and $|\cdot|$ denotes the Euclidean norm on both $\mathbb{R}^{3\times3}\cong\mathbb{R}^9$ and~$\mathbb{R}^3$. One can check that \Cref{assumption1} is satisfied for $g_\ast=(\boldsymbol{R}_\ast,\boldsymbol{r}_\ast)$ with $\boldsymbol{R}_\ast=\boldsymbol{I}_3$ and $\boldsymbol{r}_\ast=0$. The component functions of $\boldsymbol{u}=(\boldsymbol{u}_{\boldsymbol{\Omega}},\boldsymbol{u}_{\boldsymbol{V}})$ are defined by $u^i(\tau):=\sqrt{2}(7-i)\cos((7-i)\tau)$ and the function~$\alpha$ is defined by~\cref{eq:15}. We set $\omega:=5$, $h:=1$,%
\begin{equation}\label{eq:40}
\boldsymbol{R}_0 := \tfrac{1}{3}\left[\begin{smallmatrix} -1 & 2 & -2 \\ -2 & 1 & 2 \\ 2 & 2 & 1 \end{smallmatrix}\right] \in \mathrm{SO}(3), \ \ \boldsymbol{r}_0 := \left[\begin{smallmatrix} 1 \\ 1 \\ 1 \end{smallmatrix}\right] \in \mathbb{R}^3.
\end{equation}%
\Cref{figure2} shows the components $\boldsymbol{R}$ and $\boldsymbol{r}$ of the solution of~\cref{eq:36} with initial state $((\boldsymbol{R}_0,\boldsymbol{r}_0),(\boldsymbol{0},\boldsymbol{0}),(\boldsymbol{0},\boldsymbol{0}),0)$ at initial time $0$ for two different choices of the parameters $a$, $\lambda$, $b$, $\kappa$ satisfying inequality~\cref{eq:11}. We can see that $\boldsymbol{R}$ and $\boldsymbol{r}$ converge into neighborhoods of $\boldsymbol{R}_\ast$ and $\boldsymbol{r}_\ast$, respectively. The size of those neighborhoods can be made arbitrary small by choosing the control parameter $\omega$ sufficiently large. As in \Cref{figure1}, we can also observe in \Cref{figure2} that the choice of the parameters $a$, $\lambda$, $b$, $\kappa$ has a significant impact on the performance of the closed-loop system.

\section{Discussion}\label{sec:05}
An implementation of the proposed method is very energy consuming since it requires periodic forces and torques with sufficiently large amplitudes and frequencies. This can lead to practical issues in applications to systems with a limited amount of available energy, such as satellites in space. On the other hand, the term ``large-amplitude high-frequency perturbations'' should be understood in a relative sense. The employed perturbation signals may have small amplitudes and frequencies if the gains $a$, $\lambda$, $b$, $\kappa$ are sufficiently small so that the time-evolution of the state variables is much slower than the variations of the perturbation signals. This statement is confirmed by the simulation results in \Cref{sec:04}. In \Cref{sec:04:01}, the gains $a$, $\lambda$, $b$, $\kappa$ are large and we have to choose a large value for the amplitude and frequency parameter $\omega$ to ensure a good performance. In this case, we can observe fast convergence to the desired configuration. In \Cref{sec:04:02}, the gains $a$, $\lambda$, $b$, $\kappa$ are small and a moderate value of $\omega$ is already sufficient to obtain convergence into a small neighborhood of the desired configuration. However, the speed of convergence towards the desired configuration in \Cref{sec:04:02} is much slower than in \Cref{sec:04:01}. Consequently, the perturbation-based method has to be applied for a longer period of time, which means more consumption of energy.

We emphasize that our analysis is restricted to fully actuated systems, but many real-world vehicles (such as aircraft and underwater vehicles) are underactuated and they are subject to nonholonomic velocity constraints. An analysis of the method for a certain class of underactuated systems is currently in preparation. We also note that in most environments (except for outer space) damping naturally occurs. In the presence of sufficiently strong dissipation, the simple method in \cite{Suttner20221} is probably a better choice. However, if naturally concurring dissipation is rather weak, then the method in the present paper with artificial damping leads to faster convergence to a desired configuration.

%-----------------------------------------------------------------------------------------------------------------------------------------------------------------
%-----------------------------------------------------------------------------------------------------------------------------------------------------------------
%                                                                Bibliography
%-----------------------------------------------------------------------------------------------------------------------------------------------------------------
%-----------------------------------------------------------------------------------------------------------------------------------------------------------------

\bibliographystyle{plain}
\bibliography{bibFile}
\appendix

%-----------------------------------------------------------------------------------------------------------------------------------------------------------------
%-----------------------------------------------------------------------------------------------------------------------------------------------------------------
%                                                                Appendix A
%-----------------------------------------------------------------------------------------------------------------------------------------------------------------
%-----------------------------------------------------------------------------------------------------------------------------------------------------------------

\section{Proof of \texorpdfstring{\Cref{proposition1}}{Proposition~\ref{proposition1}}}\label{sec:appendixA}
We follow the proof of Theorem~1 in \cite{Moreau2000}.

Assume that the solutions of~\cref{eq:01} in the variables~\cref{eq:02} approximate the solutions of~\cref{eq:03}. Let $x_\ast$ be a point of $M$ and let $S$ be a neighborhood of $x_\ast$ in $M$. For every $r>0$, let $K_\ast(r)$ denote set of $x\in{M}$ with $|x-x_\ast|\leq{r}$. Suppose that $M$ is a topologically closed subset of the ambient Euclidean space. Then $K_\ast(r)$ is compact for every $r>0$. Assume that $x_\ast$ is $S$-uniformly asymptotically stable for~\cref{eq:03}. We prove that $x_\ast$ is $S$-practically uniformly asymptotically stable for~\cref{eq:01} in the variables~\cref{eq:02}.

\emph{Practical uniform stability.} Fix an arbitrary small $\varepsilon>0$. Since $x_\ast$ is assumed to be uniformly stable for~\cref{eq:03}, there exists some sufficiently small $\delta\in(0\vphantom{)},\varepsilon/2\vphantom{[}]$ such that every maximal solution of~\cref{eq:03} that starts in $K_\ast(\delta)$ stays in $K_\ast(\varepsilon/2)$ at any later time. After possibly shrinking $\delta>0$, we may suppose that $K_\ast(\delta)$ is contained in $S$. Set $\nu:=\delta/2$. Since $x_\ast$ is assumed to be $S$-uniformly attractive for~\cref{eq:03}, there exists some sufficiently large $\Delta>0$ such that every maximal solution of~\cref{eq:03} that starts in $K_\ast(\delta)$ enters $K_\ast(\nu)$ at latest after the time span $\Delta$ and then stays in $K_\ast(\nu)$ at any later time. Since $K_\ast(\varepsilon/2)$ is a compact subset of $M$ and since the solutions of \cref{eq:01} in the variables \cref{eq:02} are assumed to approximate the solutions of \cref{eq:03}, there exists some sufficiently large $\omega_0>0$ such that, for every $\omega\geq\omega_0$, every $t_0\in\mathbb{R}$, and every $\bar{x}_0\in{K_\ast(\delta)}$, the maximal solution $\bar{x}$ of~\cref{eq:03} with initial condition $\bar{x}(t_0)=\bar{x}_0$ and the maximal solution $x$ of~\cref{eq:01} with initial condition $x(t_0)=(\phi^{\omega}_{t_0})^{-1}(\bar{x}_0)$ satisfy $|\phi^{\omega}_t(x(t))-\bar{x}(t)|\leq\nu$ for every $t\in[t_0,t_0+\Delta]$. For the rest of this paragraph, fix  an arbitrary $\omega\geq\omega_0$, an arbitrary $t_0\in\mathbb{R}$, and an arbitrary $\tilde{x}_0\in{K_\ast(\delta)}$. For every integer $k$, set $t_k:=t_0+k\Delta$. Let $x$ be the maximal solution of~\cref{eq:01} with initial condition $x(t_0)=(\phi^{\omega}_{t_0})^{-1}(\tilde{x}_0)$. We proceed by induction on the nonnegative integer $k$. Suppose that we have already shown that $\phi^{\omega}_t(x(t))\in{K_\ast(\varepsilon)}$ for every $t\in[t_0,t_k]$ and that $\phi^{\omega}_{t_k}(x(t_k))\in{K_\ast(\delta)}$. This is trivially true for $k=0$. Now, for the induction step, set $\tilde{x}_k:=\phi^{\omega}_{t_k}(x(t_k))$. Let $\bar{x}$ be the maximal solution of~\cref{eq:03} with initial condition $\bar{x}(t_k)=\tilde{x}_k$. Then, we know that $\bar{x}(t)\in{K_\ast(\varepsilon/2)}$ for every $t\in[t_k,t_{k+1}]$ and that $\bar{x}(t_{k+1})\in{K_\ast(\delta/2)}$. By the triangle inequality, it follows that $\phi^{\omega}_t(x(t))\in{K_\ast(\varepsilon)}$ for every $t\in[t_k,t_{k+1}]$ and that $\phi^{\omega}_{t_{k+1}}(x(t_{k+1}))\in{K_\ast(\delta)}$. The proof of practical uniform stability is complete.
	
\emph{Practical uniform attraction.} Fix an arbitrary small $\varepsilon>0$ and an arbitrary large $r>0$. After possibly shrinking $\varepsilon>0$, we may suppose that $\varepsilon\leq{r}$. We already know from the preceding paragraph that $x_\ast$ is practically uniformly stable for~\cref{eq:01} in the variables~\cref{eq:02}. Therefore, there exists some sufficiently small $\delta>0$ and some sufficiently large $\omega_1>0$ such that, for every $\omega\geq\omega_1$, every $t_0\in\mathbb{R}$, and every $\tilde{x}_0\in{K_\ast(\delta)}$, the maximal solution $x$ of~\cref{eq:01} with initial condition $x(t_0)=(\phi^{\omega}_{t_0})^{-1}(\tilde{x}_0)$ satisfies $\phi^{\omega}_t(x(t))\in{K_\ast(\varepsilon)}$ for every $t\geq{t_0}$. Since $x_\ast$ is assumed to be $S$-uniformly attractive for~\cref{eq:03}, there exist sufficiently large $\Delta,\bar{R}>0$ such that, for every $t_0\in\mathbb{R}$ and every $\bar{x}_0\in{S}$ with $\bar{x}_0\in{K_\ast(r)}$, the maximal solution $\bar{x}$ of~\cref{eq:03} with initial condition $\bar{x}(t_0)=\bar{x}_0$ satisfies $\bar{x}(t)\in{K_\ast(\bar{R})}$ for every $t\geq{t_0}$ and $\bar{x}(t)\in{K_\ast(\delta/2)}$ for every $t\geq{t_0+\Delta}$. Set $\nu:=\delta/2$. Since $K_\ast(\bar{R})$ is a compact subset of $M$ and since the solutions of \cref{eq:01} in the variables \cref{eq:02} are assumed to approximate the solutions of \cref{eq:03}, there exists some sufficiently large $\omega_2>0$ such that, for every $\omega\geq\omega_2$, every $t_0\in\mathbb{R}$, and every $\bar{x}_0\in{K_\ast(r)}$, the maximal solution $\bar{x}$ of~\cref{eq:03} with initial condition $\bar{x}(t_0)=\bar{x}_0$, and the maximal solution $x$ of~\cref{eq:01} with initial condition $x(t_0)=(\phi^{\omega}_{t_0})^{-1}(\bar{x}_0)$ satisfy $|\phi^{\omega}_t(x(t))-\bar{x}(t)|\leq\nu$ for every $t\in[t_0,t_0+\Delta]$. Let $\omega_0>0$ be the maximum of $\omega_1$ and $\omega_2$, and set $R:=\bar{R}+\nu$. For the rest of this paragraph, fix an arbitrary $\omega\geq\omega_0$, an arbitrary $t_0\in\mathbb{R}$, and an arbitrary $\tilde{x}_0\in{K_\ast(r)}$. Set $t_1:=t_0+\Delta$. Let $x$ be the maximal solution of~\cref{eq:01} with initial condition $x(t_0)=(\phi^{\omega}_{t_0})^{-1}(\tilde{x}_0)$ and let $\bar{x}$ be the maximal solution of~\cref{eq:03} with initial condition $\bar{x}(t_0)=\tilde{x}_0$. Then, we know that $\bar{x}(t)\in{K_\ast(\bar{R})}$ for every $t\in[t_0,t_1]$ and that $\bar{x}(t_1)\in{K_\ast(\delta/2)}$. By the triangle inequality, it follows that $\phi^{\omega}_t(x(t))\in{K_\ast(R)}$ for every $t\in[t_0,t_1]$ and that $\phi^{\omega}_{t_1}(x(t_1))\in{K_\ast(\delta)}$. Consequently, $\phi^{\omega}_t(x(t))\in{K_\ast(\varepsilon)}$ and in particular $\phi^{\omega}_t(x(t))\in{K_\ast(R)}$ for every $t\geq{t_1}$. The proof of practical uniform attraction is complete.

%-----------------------------------------------------------------------------------------------------------------------------------------------------------------
%-----------------------------------------------------------------------------------------------------------------------------------------------------------------
%                                                                Appendix B
%-----------------------------------------------------------------------------------------------------------------------------------------------------------------
%-----------------------------------------------------------------------------------------------------------------------------------------------------------------

\section{Proof of \texorpdfstring{\Cref{proposition2}}{Proposition~\ref{proposition2}}}\label{sec:appendixB}
We prove a slightly more general approximation result than \Cref{proposition2}; see \Cref{eq:appendix2Lemma} below. The general approximation result contains \Cref{proposition2} as a special case; see \Cref{eq:appendix2Remark} below.

Let $Q$ be an embedded submanifold of Euclidean space and let $\Xi$ be a vector subspace of Euclidean space. Let $TQ$ denote the tangent bundle of $Q$ and let $M:=Q\times\Xi$ be the product manifold of $Q$ and $\Xi$, which is again an embedded submanifold of Euclidean space. Let $B\colon\Xi\times\Xi\to\Xi$ be a bilinear map. For every differentiable curve $\xi\colon{I}\to\Xi$, we define $\nabla_\xi\xi\colon{I}\to\Xi$ by%
\begin{equation}\label{eq:B.01}
\nabla_\xi\xi(t) \ := \ \dot{\xi}(t) + B(\xi(t),\xi(t)).
\end{equation}%
For every point $q$ of $Q$, let $L(q)$ be a linear map from $\Xi$ into the tangent space $T_qQ$ to $Q$ at $q$. We suppose that the map $M\to{TQ}$, $(q,\xi)\mapsto{L(q)\xi}$ is smooth. For all smooth maps $Y,Z\colon{Q}\to\Xi$, we define smooth maps $\nabla_YZ,\langle{Y\colon\!Z}\rangle\colon{Q}\to\Xi$ by%
\begin{align}
\nabla_{Y}Z(q) & \ := \ \mathrm{d}Z(q)L(q)Y(q) + B(Y(q),Z(q)), \label{eq:B.02} \allowdisplaybreaks \\
\langle{Y\colon\!Z}\rangle(q) & \ := \ \nabla_{Y}Z(q) + \nabla_{Z}Y(q), \label{eq:B.03}
\end{align}%
where $\mathrm{d}Z(q)\colon{T_qQ}\to\Xi$ denotes the differential of $Z$ at $q$. For every point $q$ of $Q$, let $R(q)$ be an endomorphism of $\Xi$. We suppose that the map $M\to\Xi$, $(q,\xi)\mapsto{R(q)\xi}$ is smooth. Let $X\colon{M}\to{TQ}$ be a smooth map such that, for every $\xi\in\Xi$, the map $X(\cdot,\xi)\colon{Q}\to{TQ}$ is a vector field on $Q$. Let $Y_0\colon{M}\to\Xi$ and $Y_1,\ldots,Y_m\colon{Q}\to\Xi$ be smooth maps. Let $T$ be a positive real number and let $u^1,\ldots,u^m$ be measurable and bounded real-valued functions on $\mathbb{R}$. We suppose that $u^1,\ldots,u^m$ are $T$-periodic and zero-mean. For each $i\in\{1,\ldots,m\}$, let $U^i$ be the zero-man antiderivative of $u^i$. For all $i,j\in\{1,\ldots,m\}$, we define the real number%
\begin{equation}\label{eq:B.04}
\Lambda^{ij} \ := \ \frac{1}{2T}\int_0^TU^i(\tau)\,U^j(\tau)\,\mathrm{d}\tau.
\end{equation}%
We will see that the solutions of the system%
\begin{subequations}\label{eq:B.05}%
\begin{align}
\dot{q} & \ = \ X\Big(q,\xi - \sum_{i=1}^mU^i(\omega{t})\,Y_i(q)\Big) + L(q)\xi, \allowdisplaybreaks \\
\nabla_\xi\xi & \ = \ Y_0\Big(q,\xi - \sum_{i=1}^mU^i(\omega{t})\,Y_i(q)\Big) + R(q)\xi \allowdisplaybreaks \\
& \qquad + \sum_{i=1}^m\omega\,u^i(\omega{t})\,Y_i(q)
\end{align}%
\end{subequations}%
on $M=Q\times\Xi$ in the variables%
\begin{subequations}\label{eq:B.06}%
\begin{align}
\tilde{q} & \ = \ q, \allowdisplaybreaks \\
\tilde{\xi} & \ = \ \xi - \sum_{i=1}^mU^i(\omega{t})\,Y_i(q)
\end{align}%
\end{subequations}%
approximate the solutions of the system%
\begin{subequations}\label{eq:B.07}%
\begin{align}
\dot{\bar{q}} & \ = \ X(\bar{q},\bar{\xi}) + L(\bar{q})\bar{\xi}, \allowdisplaybreaks \\
\nabla_{\bar{\xi}}\bar{\xi} & \ = \ Y_0(\bar{q},\bar{\xi}) + R(\bar{q})\bar{\xi} - \sum_{i=1}^m\Lambda^{ij}\langle{Y_i\colon\!Y_j}\rangle(\bar{q})
\end{align}%
\end{subequations}%
on $M$ in the sense of \Cref{definition4}.%
\begin{lemma}\label{eq:appendix2Lemma}
The solutions of \cref{eq:B.05} in the variables \cref{eq:B.06} approximate the solutions of \cref{eq:B.07}.
\end{lemma}%
\begin{proof}
We follow the proof of Theorem~4.1 in~\cite{Bullo2002} and the proof of Theorem~9.32 in~\cite{BulloBook}.

For all $i,j\in\{1,\ldots,m\}$, define $v^i,v^{ij}\colon\mathbb{R}\to\mathbb{R}$ by%
\begin{align*}
v^i(\tau) & \ := \ U^i(\tau), \allowdisplaybreaks \\
v^{ij}(\tau) & \ := \ \tfrac{1}{2}\,U^i(\tau)\,U^j(\tau) - \Lambda^{ij},
\end{align*}%
define $f_0^1,f_i^1,f_{ij}^1\colon{M}\to{TQ}$ by%
\begin{align*}
f_0^1(\tilde{q},\tilde{\xi}) & \ := \ X(\tilde{q},\tilde{\xi}) + L(\tilde{q})\tilde{\xi}, \allowdisplaybreaks \\
f_i^1(\tilde{q},\tilde{\xi}) & \ := \ L(\tilde{\xi})Y_i(\tilde{\xi}), \allowdisplaybreaks \\
f_{ij}^1(\tilde{q},\tilde{\xi}) & \ := \ 0,
\end{align*}%
and define $f_0^2,f_i^2,f_{ij}^2\colon{M}\to\Xi$ by%
\begin{align*}
f_0^2(\tilde{q},\tilde{\xi}) & \ := \ Y_0(\tilde{q},\tilde{\xi}) + R(\tilde{q})\tilde{\xi} - \sum_{i,j=1}^m\Lambda^{ij}\,\langle{Y_i\colon\!Y_j}\rangle(\tilde{q}), \allowdisplaybreaks \\
f_i^2(\tilde{q},\tilde{\xi}) & \ := \ R(\tilde{q})Y_i(\tilde{q}) - B(\tilde{\xi},Y_i(\tilde{q})) - B(Y_i(\tilde{q}),\tilde{\xi}) \allowdisplaybreaks \\
& \qquad - \mathrm{d}Y_i(\tilde{q})X(\tilde{q},\tilde{x}) - \mathrm{d}Y_i(\tilde{q})L(\tilde{q})\tilde{\xi}, \allowdisplaybreaks \\
f_{ij}^2(\tilde{q},\tilde{\xi}) & \ := \ -\langle{Y_i\colon\!Y_j}\rangle(\tilde{q}).
\end{align*}%
Let $J$ be the set of indices $i$ and $ij$ with $i,j\in\{1,\ldots,m\}$. For each $\iota\in\{0\}\cup{J}$, let $f_\iota$ be the vector field on $M$ with components $f_\iota^1$ and $f_\iota^2$. A direct computation shows that \cref{eq:B.05} in the variables \cref{eq:B.06} is the system%
\begin{equation}\label{eq:B.08}
\dot{\tilde{x}} \ = \ f_0(\tilde{x}) + \sum_{\iota\in{J}}v^\iota(\omega{t})\,f_\iota(\tilde{x})
\end{equation}%
on $M$. The proof is complete if we can show that the solutions of \cref{eq:B.08} approximate the solutions of%
\begin{equation}\label{eq:B.09}
\dot{\bar{x}} \ = \ f_0(\bar{x}),
\end{equation}%
since \cref{eq:B.09} is the same as \cref{eq:B.07}.

For each $\iota\in{J}$, the function $v^\iota$ is $T$-periodic and zero-mean, and therefore there exists a zero-mean antiderivative $V^\iota$ of $v^\iota$. Let $\varphi$ be the smooth inclusion map from $M$ into the ambient Euclidean space. For every vector field $f$ on $M$, let $f\varphi$ denote the component-wise Lie derivative of $\varphi$ with respect to $f$. Define time-dependent vector-valued functions $D_1\varphi$ and $D_2\varphi$ on $M$ by%
\begin{align*}
(D_1\varphi)(\tau,\tilde{x}) & \ := \ \sum_{\iota\in{J}}V^\iota(\tau)\,(f_\iota\varphi)(\tilde{x}), \allowdisplaybreaks \\
(D_2\varphi)(\tau,\tilde{x}) & \ := \ \sum_{\iota\in{J}}V^\iota(\tau)\,(f_0(f_\iota\varphi))(\tilde{x}) \allowdisplaybreaks \\
& \qquad + \sum_{\nu,\iota\in{J}}v^\nu(\tau)\,V^\iota(\tau)\,(f_\nu(f_\iota\varphi))(\tilde{x}).
\end{align*}%
For every vector-valued map $a$ on an interval $I$ and for all $t_1,t_2\in{I}$, we use the notation $[a(t)]_{t=t_1}^{t=t_2}$ for the difference $a(t_2)-a(t_1)$. A direct computation, using integration by parts, shows that every solution $\tilde{x}\colon{I}\to{M}$ of \cref{eq:B.08} satisfies the integral equation%
\begin{align*}
& \varphi(\tilde{x}(t_2)) \ = \ \varphi(\tilde{x}(t_1)) + \int_{t_1}^{t_2}(f_0\varphi)(\tilde{x}(t))\,\mathrm{d}t \\
& + \frac{1}{\omega}\Big[(D_1\varphi)(\omega{t},\tilde{x}(t))\Big]_{t=t_1}^{t=t_2} - \frac{1}{\omega}\int_{t_1}^{t_2}(D_2\varphi)(\omega{t},\tilde{x}(t))\,\mathrm{d}t
\end{align*}%
for all $t_1,t_2\in{I}$.

To prove the asserted approximation property, fix an arbitrary large compact subset $K$ of $M$, an arbitrary small approximation error $\nu>0$ and an arbitrary large time span $\Delta>0$. After possibly reducing $\nu>0$, we can ensure that the closed $\nu$-neighborhood $\tilde{K}$ of $K$ in $M$ is compact. Since $f_0\varphi$ is smooth, there exists some sufficiently large Lipschitz constant $L>0$ such that
\[
\big|(f_0\varphi)(\tilde{x}) - (f_0\varphi)(\bar{x})\big| \ \leq \ L\,|\varphi(\tilde{x}) - \varphi(\bar{x})|
\]
for all $\tilde{x},\bar{x}\in\tilde{K}$. Since $D_1\varphi$ and $D_2\varphi$ are $T$-periodic in the first argument, there exist sufficiently large bounds $c_1>0$ and $c_2>0$ such that
\[
|(D_1\varphi)(\tau,\tilde{x})| \ \leq \ c_1 \qquad \text{and} \qquad |(D_2\varphi)(\tau,\tilde{x})| \ \leq \ c_2
\]
for every $\tau\in\mathbb{R}$ and every $\tilde{x}\in{K}$. Choose some
\[
\omega_0 \ > \ \tfrac{1}{\nu}\,(2\,c_1+\Delta\,c_2)\,\mathrm{e}^{\Delta{L}}.
\]
We show that this choice of $\omega_0$ does the job. To this end, fix an arbitrary initial time $t_0\in\mathbb{R}$ and an arbitrary initial value $\bar{x}_0$. Suppose that the maximal solution $\bar{x}$ of \cref{eq:B.09} with initial condition $\bar{x}(t_0)=\bar{x}_0$ satisfies $\bar{x}(t)\in{K}$ for every $t\in[t_0,t_0+\Delta]$. Fix an arbitrary $\omega\geq\omega_0$ and let $\tilde{x}\colon\tilde{I}\to{M}$ be the maximal solution of \cref{eq:B.08} with initial condition $\tilde{x}(t_0)=\bar{x}_0$. Let $t_1$ be the largest element in the intersection of $\tilde{I}$ and $[t_0,t_0+\Delta]$ such that $\tilde{x}(t)\in\tilde{K}$ for every $t\in[t_0,t_1]$. Using the integral equation for the solutions of \cref{eq:B.08} and the Gronwall lemma in integral form, we obtain that
\[
\big|\varphi(\tilde{x}(t)) - \varphi(\bar{x}(t))\big| \ \leq \ \tfrac{1}{\omega}\,(2\,c_1+\Delta\,c_2)\,\mathrm{e}^{\Delta{L}} \ < \ \nu
\]
for every $t\in[t_0,t_1]$. In particular, this implies $t_1=t_0+\Delta$ and completes the proof.%
\end{proof}%
\begin{remark}\label{eq:appendix2Remark}\em
To see that \Cref{proposition2} is indeed a special case of \Cref{eq:appendix2Lemma}, set $Q:=G\times\mathbb{R}$ and $\Xi:=\mathfrak{g}\times\mathbb{R}^m$. The maps $B$, $L$, $R$, $X$, $Y_0$, and $Y_i$ with $i\in\{1,\ldots,m\}$ are defined by%
\begin{align*}
B((v,\boldsymbol{w}),(v',\boldsymbol{w}')) & \ := \ \big(\stackrel{\mathfrak{g}}{\nabla}_v\!v',\boldsymbol{0}\big), \allowdisplaybreaks \\
L(g,\eta)(v,r) & \ := \ \big(T_eL_g(v),0\big), \allowdisplaybreaks \\
R(g,\eta)(v,\boldsymbol{w}) & \ := \ \big(-Rv,\boldsymbol{0}\big), \allowdisplaybreaks \\
X((g,\eta),(\tilde{v},\tilde{\boldsymbol{w}})) & \ := \ \big(0,\,-h\,\eta+h\,\psi(g)\big), \allowdisplaybreaks \\
Y_0((g,\eta),(\tilde{v},\tilde{\boldsymbol{w}})) & \ := \ \big(\lambda^2\,\alpha^2(\psi(g)-\eta)\,\mu,\boldsymbol{0}\big) \allowdisplaybreaks \\
& \qquad - \big(b\,\tilde{w},\,a\,\tilde{\boldsymbol{w}}\big), \allowdisplaybreaks \\
Y_i(g,\eta) & \ := \ \alpha(\psi(g)-\eta)\,\big(\lambda\,e_i, \kappa\,\boldsymbol{e}_i\big),
\end{align*}%
where $m:=n$ is the dimension of $G$. Equation~\cref{eq:13} implies that $\Lambda_{ij}=1/2$ for $i=j$ and $\Lambda_{ij}=0$ otherwise. Now it is easy to verify that \cref{eq:20}, \cref{eq:22}, and \cref{eq:23} are special cases of \cref{eq:B.05}, \cref{eq:B.06}, and \cref{eq:B.07}.
\end{remark}%

%-----------------------------------------------------------------------------------------------------------------------------------------------------------------
%-----------------------------------------------------------------------------------------------------------------------------------------------------------------
%                                                                Appendix C
%-----------------------------------------------------------------------------------------------------------------------------------------------------------------
%-----------------------------------------------------------------------------------------------------------------------------------------------------------------

\section{Proofs of \texorpdfstring{\Cref{theorem1,theorem2}}{Theorems~\ref{theorem1} and~\ref{theorem2}}}\label{sec:appendixC}
Recall that we use \Cref{notation1,notation2} and that the real-valued function $V$ on $M=G\times\mathfrak{g}\times\mathbb{R}^n\times\mathbb{R}$ is defined by \cref{eq:28}. Let $\dot{V}$ denote the derivative of $V$ along solutions of~\cref{eq:23}. A direct computation shows that $\dot{V}$ is given by%
\begin{subequations}\label{eq:C.01}%
\begin{align}
\dot{V}(g,v,& \boldsymbol{w},\eta) \ = \ -h\,\lambda^2\,\beta'(\psi(g)-\eta)\,(\psi(g)-\eta) \label{eq:C.01:a} \allowdisplaybreaks \\
& - \|v\|_R^2 - \|v\|_{R\Xi}^2 - \tfrac{\lambda}{\kappa^2}\,(a\lambda-b\kappa)\,\|w\|_\Xi^2 \label{eq:C.01:b} \allowdisplaybreaks \\
& - \mathbb{I}\big(\Xi\stackrel{\mathfrak{g}}{\nabla}_v\!v,v-\tfrac{\lambda}{\kappa}\,w\big). \label{eq:C.01:c}
\end{align}%
\end{subequations}%

%-----------------------------------------------------------------------------------------------------------------------------------------------------------------
%                                                                Appendix C.1
%-----------------------------------------------------------------------------------------------------------------------------------------------------------------

\subsection{Proof of \texorpdfstring{\Cref{theorem1}}{Theorem~\ref{theorem1}}}
Suppose that \Cref{assumption1} is satisfied with $g_\ast$ as therein. For every $\varepsilon>0$, define a real-valued function $V_\varepsilon$ on $M$ by%
\begin{align*}
V_{\varepsilon}(g,v,\boldsymbol{w},\eta) & \ := \ V(g,v,\boldsymbol{w},\eta) - \varepsilon\,\tfrac{b\kappa}{\lambda}\,\mathbb{I}(\mathfrak{grad}\psi(g),v) \allowdisplaybreaks \\
& \qquad + \varepsilon\,\tfrac{a\lambda}{\kappa}\,\mathbb{I}(\mathfrak{grad}\psi(g),w).
\end{align*}%
Such a modification of the total energy is sometimes referred to as \emph{Chetaev's trick} (cf.~Remark 6.46-3 in~\cite{BulloBook}).

Set $x_\ast:=(g_\ast,0,\boldsymbol{0},\psi(g_\ast))\in{M}$. Note that $V_\varepsilon(x_\ast)=0$ for every $\varepsilon>0$. We show that, for every sufficiently small $\varepsilon>0$, the function $V_\varepsilon$ is locally positive-definite about $x_\ast$. Using the Cauchy-Schwarz inequality and completion of squares, it is easy to see that there exist $\vartheta_1,\vartheta_2>0$ such that the kinetic energy terms of $V$ satisfy the estimate%
\begin{equation*}
\tfrac{1}{2}\,\|v\|^2 + \tfrac{1}{2}\,\|v-\tfrac{\lambda}{\kappa}\,w\|_\Xi^2 \ \geq \ \vartheta_1\,\|v\|^2 + \vartheta_2\,\|w\|^2
\end{equation*}%
for all $v,w\in\mathfrak{g}$. It follows from \Cref{assumption1} and Proposition~6.30 in~\cite{BulloBook} that there exist $\gamma>0$ and some sufficiently small neighborhood $U$ of $g_\ast$ in $G$ such that%
\begin{equation*}
\psi(g) - \psi(g_\ast) \ \geq \ \gamma\,\|\mathfrak{grad}\psi(g)\|^2 \ > \ 0
\end{equation*}%
for every $g\in{U}$ with $g\neq{g_\ast}$. Using the Cauchy-Schwarz inequality and completion of squares, we obtain that%
\begin{align*}
& V_\varepsilon(g,v,\boldsymbol{w},\eta) \ \geq \ \lambda^2\,\beta(\psi(g)-\eta) \allowdisplaybreaks \\
& \qquad + \big(\lambda^2\,(\alpha\alpha')(0)\,\gamma - \varepsilon\,\tfrac{b\kappa}{\lambda} - \varepsilon\,\tfrac{a\lambda}{\kappa}\big)\,\|\mathfrak{grad}\psi(g)\|^2 \allowdisplaybreaks \\
& \qquad + \big(\vartheta_1 - \varepsilon\,\tfrac{b\kappa}{\lambda}\big)\|v\|^2 + \big(\vartheta_2 - \varepsilon\,\tfrac{a\lambda}{\kappa}\big)\|w\|^2
\end{align*}%
for every $\varepsilon>0$ and every $(g,v,\boldsymbol{w},\eta)\in{M}$ with $g\in{U}$. Now, using that $\beta$ is positive definite about $0$, one can show that there exists some sufficiently small $\varepsilon_1>0$ such that $V_\varepsilon(x)>0$ for every $\varepsilon\in(0,\varepsilon_1)$ and every $x=(g,v,\boldsymbol{w},\eta)\in{M}$ with $g\in{U}$ and $x\neq{x_\ast}$.

For every $g\in{G}$, define the endomorphism
\[
\mathfrak{Hess}\psi(g) \ := \ (T_eL_g)^{-1}\circ\mathrm{Hess}\psi(g)\circ{L_{g\ast{e}}}
\]
of $\mathfrak{g}$, where $\mathrm{Hess}\psi$ is the Hessian operator of $\psi$ with respect to the left-invariant Riemannian metric on $G$ induced by $\mathbb{I}$. For every $\varepsilon>0$, let $\dot{V}_\varepsilon$ denote the derivative of $V_\varepsilon$ along solutions of~\cref{eq:23}. A direct computation shows that $\dot{V}_\varepsilon/\varepsilon$ is given by%
\begin{align*}
& \tfrac{1}{\varepsilon}\,\dot{V}_\varepsilon(g,v,\boldsymbol{w},\eta) \ = \ \tfrac{1}{\varepsilon}\,\dot{V}(g,v,\boldsymbol{w},\eta) - \tfrac{b\kappa}{\lambda}\,\|v\|_{\mathfrak{Hess}\psi(g)}^2 \allowdisplaybreaks \\
& \quad - \lambda\,(a\,\lambda-b\,\kappa)\,(\alpha\alpha')\big(\psi(g)-\eta\big)\,\|\mathfrak{grad}\psi(g)\|^2 \allowdisplaybreaks \\
& \quad + \tfrac{a\lambda}{\kappa}\,\mathbb{I}(\mathfrak{Hess}\psi(g)v,w) - \tfrac{(a\lambda)^2-(b\kappa)^2}{\lambda\,\kappa}\,\mathbb{I}(\mathfrak{grad}\psi(g),w) \allowdisplaybreaks \\
& \quad + \tfrac{a\lambda}{\kappa}\,\mathbb{I}\big(\mathfrak{grad}\psi(g),\stackrel{\mathfrak{g}}{\nabla}_v\!w\big) + \tfrac{b\kappa}{\lambda}\,\mathbb{I}(\mathfrak{grad}\psi(g),Rv).
\end{align*}%
Note that $\dot{V}_\varepsilon(x_\ast)=0$ for every $\varepsilon>0$. We show that, for every sufficiently small $\varepsilon>0$, the function $\dot{V}_\varepsilon$ is locally negative-definite about $x_\ast$. Because of \Cref{assumption1}, after possibly shrinking the neighborhood $U$ of $g_\ast$ from the previous paragraph, we can find $h_1,h_2>0$ such that%
\begin{equation*}
\|v\|_{\mathfrak{Hess}\psi(g)}^2 \geq h_1\,\|v\|^2, \ \big|\mathbb{I}(\mathfrak{Hess}\psi(g)v,w)\big| \leq h_2\,\|v\|\,\|w\|
\end{equation*}%
for every $g\in{U}$ and all $v,w\in\mathfrak{g}$. Since $R$ is an endomorphism of $\mathfrak{g}$ and since the restriction of $\nabla$ to $\mathfrak{g}$ is bilinear, there exist $r,s>0$ such that%
\begin{equation*}
\big|\mathbb{I}(R{v},w)\big| \leq r\,\|v\|\,\|w\|, \qquad \big\|\stackrel{\mathfrak{g}}{\nabla}_v\!w\big\| \ \leq \ s\,\|v\|\,\|w\|
\end{equation*}%
for all $v,w\in\mathfrak{g}$. Since $\Xi$ is symmetric and positive-definite with respect to $\mathbb{I}$, there exist $\xi_1,\xi_2>0$ such that%
\begin{equation*}
\|v\|_\Xi^2 \ \geq \ \xi_1\,\|v\|^2, \qquad \big|\mathbb{I}(\Xi{v},w)\big| \ \leq \ \xi_2\,\|v\|\,\|w\|
\end{equation*}%
for all $v,w\in\mathfrak{g}$. Using the Cauchy-Schwarz inequality and completion of squares, one can show that there exist positive constants $c_0,\ldots,c_5$ such that%
\begin{align*}
& \dot{V}_\varepsilon(g,v,\boldsymbol{w},\eta) \allowdisplaybreaks \\
& \ \leq \ -h\,\lambda^2\,\beta'(\psi(g)-\eta)\,(\psi(g)-\eta) - c_0\,\|w\|^2 \allowdisplaybreaks \\
& \qquad - \varepsilon\,\big(c_1-\varepsilon\,c_2-\varepsilon\,c_3\,\|\mathfrak{grad}\psi(g)\|^2 \allowdisplaybreaks \\
& \qquad\qquad\qquad\qquad - \tfrac{1}{\varepsilon}\,\xi_2\,s\,\|v-\tfrac{\lambda}{\kappa}\,w\|\big)\|v\|^2 \allowdisplaybreaks \\
& \qquad - \varepsilon\,\big(c_4\,(\alpha\alpha')(\psi(g)-\eta)-\varepsilon\,c_5\big)\,\|\mathfrak{grad}\psi(g)\|^2
\end{align*}%
for every $\varepsilon>0$ and every $(g,v,\boldsymbol{w},\eta)\in{M}$ with $g\in{U}$. Now, using that $z\mapsto\beta'(z)z$ is positive definite about $0$, one can show that there exists some sufficiently small $\varepsilon_2>0$ such that, for every $\varepsilon\in(0,\varepsilon_2)$, there exists a neighborhood $W(\varepsilon)$ of $x_\ast$ in $M$ such that $\dot{V}_\varepsilon(x)<0$ for every $x\in{W(\varepsilon)}$ with $x\neq{x_\ast}$.

In summary, we have shown that, for every sufficiently small $\varepsilon>0$, the function $V_\varepsilon$ is locally positive-definite about $x_\ast$ and $\dot{V}_\varepsilon$ is locally negative-definite about $x_\ast$. By a standard Lyapunov stability criterion (see, e.g., Theorem~6.14 in~\cite{BulloBook}), this implies that the equilibrium point $x_\ast$ is locally (uniformly) asymptotically stable for~\cref{eq:23}. Now the claim follows from \Cref{proposition1,proposition2}.

%-----------------------------------------------------------------------------------------------------------------------------------------------------------------
%                                                                Section C.2
%-----------------------------------------------------------------------------------------------------------------------------------------------------------------

\subsection{Proof of \texorpdfstring{\Cref{theorem2}}{Theorem~\ref{theorem2}}}
Suppose that \Cref{assumption2,assumption3} are satisfied with $g_\ast$ and $y_0$ as therein. Set $x_\ast:=(g_\ast,0,\boldsymbol{0},\psi(g_\ast))$ and $V_0:=\lambda^2(\alpha\alpha')(0)y_0$. The claim follows from \Cref{proposition1,proposition2} if we can show that $x_\ast$ is $V^{-1}(\leq{V_0},x_\ast)$-uniformly asymptotically stable for~\cref{eq:23}. Our goal is to apply the LaSalle invariance principle (see, e.g., Theorem~6.19 in~\cite{BulloBook}). It follows from \Cref{assumption3} that $A:=V^{-1}(\leq{V_0},x_\ast)$ is a compact subset of $M$. Because of \Cref{assumption2} and equation~\cref{eq:C.01}, we have $\dot{V}\leq0$ on $M$. This implies that the sublevel set $A$ of $V$ is positively invariant for~\cref{eq:23}. Let $C$ be the set of $x\in{A}$ with $\dot{V}(x)=0$. Let $B$ be the largest positively invariant subset of $C$ for~\cref{eq:23}. It is clear $x_\ast$ is contained in $B$. Since~\cref{eq:23} is an autonomous system, the claim that $x_\ast$ is $A$-uniformly asymptotically stable for~\cref{eq:23} follows from the LaSalle invariance principle if we can show that there is no other point than $x_\ast$ in $B$. For this purpose, let $(g,v,\boldsymbol{w},\eta)$ be a solution of~\cref{eq:23} with values in $C$. Then, we obtain from \Cref{assumption2} and equation~\cref{eq:C.01} that $\psi(g)=\eta$ and $\boldsymbol{w}\equiv0$. It follows from equation~\cref{eq:23:d} that $\mathrm{grad}\psi(g)\equiv0$. Because of \Cref{assumption3}, this in turn implies that $g\equiv{g_\ast}$, and therefore $v\equiv0$ and $\eta\equiv\psi(g_\ast)$. Thus, we have shown that $(g,v,\boldsymbol{w},\eta)$ is identically equal to $x_\ast$, which completes the proof.
\end{document}